\documentclass[12pt]{article}

\usepackage{amsmath}
\usepackage{amsfonts}
\usepackage{latexsym}
\usepackage{graphicx}

\newtheorem{proposition}{Proposition}[section]

\newtheorem{corollary}{Corollary}[section]
\newtheorem{remark}{Remark}[section]

\date{ }

\begin{document}

\title{Implicit parametrizations in shape optimization: boundary observation }

\author{Cornel Marius Murea$^1$, Dan Tiba$^2$\\
{\normalsize $^1$ D\'epartement de Math\'ematiques, IRIMAS,}\\
{\normalsize Universit\'e de Haute Alsace, France,}\\
{\normalsize cornel.murea@uha.fr}\\
{\normalsize $^2$ Institute of Mathematics (Romanian Academy) and}\\ 
{\normalsize Academy of Romanian Scientists, Bucharest, Romania,}\\ 
{\normalsize dan.tiba@imar.ro}
}

\maketitle

\begin{abstract}
We present first a brief review of the existing literature on shape optimization, stressing the recent use of Hamiltonian systems in topology optimization. In the second section, we collect some preliminaries on the implicit parametrization theorem, especially in dimension two, which is a case of interest in shape optimization. The formulation of the problem is also discussed. The approximation via penalization and its differentiability properties are analyzed in Section 3. Next, we investigate the discretization process in Section 4. The last section is devoted to numerical experiments.

\vspace{3mm}

\vspace{2mm}
\textbf{MCS 2020:} 49M20; 49Q10
\end{abstract}

\section{Introduction}
\setcounter{equation}{0}

Geometric optimization problems have a very long history
(we mention just the Dido's problem, almost three thousands years old, Kline \cite{Kline1972}), but 
shape optimization problems are a relatively young development of the calculus of variations. There exist already some very good monographs,
 Pironneau \cite{Pironneau1984},  Haslinger and Neittaanm\"aki \cite{Haslinger1996},
Sokolowski and Zolesio \cite{Sokolowski1992}, Delfour and Zolesio \cite{Delfour2001},
Neittaanm\"aki, Sprekels and Tiba \cite{NS_Tiba2006}, Bucur and Buttazzo \cite{Bucur2005},
Henrot and Pierre \cite{Henrot2005}, devoted to this subject. In general, just certain types of boundary variations for the unknown domains, are taken into account. The well known level set method, \cite{OS}, \cite{OF}, \cite{A}, \cite{MAJ}, investigates topological optimization questions as well, both from the theoretical and numerical points of view. We underline that our approach combines boundary and topological variations and is essentially different from the level set method, although level functions are used (for instance the Hamilton-Jacobi equation is not necessary here - we just use ordinary differential Hamiltonian systems, etc.).

A typical example of shape optimization problem, defined on a given family 
$\mathcal{O}$ of bounded domains $\Omega\in\mathcal{O}$, $\Omega\subset D\subset \mathbb{R}^d$,
looks as follows:
\begin{eqnarray}
  \min_{\Omega\in\mathcal{O}}
  \int_\Lambda j\left( \mathbf{x}, y_\Omega(\mathbf{x}) \right) d\mathbf{x} ,
\label{1.1}\\
-\Delta y_\Omega =f\hbox{ in }\Omega,
\label{1.2}\\
y_\Omega = 0\hbox{ on }\partial\Omega.
\label{1.3}
\end{eqnarray}
Other boundary conditions, other differential operators or cost functionals may be
as well considered in (\ref{1.1})-(\ref{1.3}). Supplementary constraints on $\Omega$
or $y_\Omega$ may be also imposed.

Above, $\Lambda$ may be $\Omega$ or some part of $\Omega$, or it may be
$\partial\Omega$ or some part of $\partial\Omega$. The functional 
$j(\cdot,\cdot):\Omega\times\mathbb{R}\rightarrow \mathbb{R}$ is Carath\'eodory, $f\in L^p(D)$, $p>2$.
The cost may also depend on $\nabla y_\Omega$ in certain situations.
Regularity assumption on $\Omega\in\mathcal{O}$, other assumptions, will be
imposed in the sequel, when necessity appears.

Shape optimization problems (\ref{1.1})-(\ref{1.3}) have a similar structure with
an optimal control problem, but the minimization parameter is the domain $\Omega$
itself, where the problem is defined.

In optimal control theory, boundary observation is an important and realistic case and
this paper is devoted to the study of boundary cost functionals in optimal design theory.
Special cases of this type have been already considered by Pironneau \cite{Pironneau1984},
Haslinger and Neittaanm\"aki \cite{Haslinger1996}, Sokolowski and Zolesio \cite{Sokolowski1992}.

The recent implicit parametrization approach, using Hamiltonian systems developped by
Tiba \cite{Tiba2013}, \cite{Tiba2018}, Nicolai and Tiba \cite{N_Tiba2015} offers a new way of handling effectively
boundary cost integrals and clarifies regularity questions, allowing developments up to
numerical experiments. Related results can be found in Tiba \cite{Tiba2018b},
\cite{Tiba2016}, \cite{MT2019a}, where the employed methodology is based on the penalization of the Dirichlet problem, but also uses the representation of
the unknown geometry via Hamiltonian systems. The family of unknown admissible domains is very general
 and the {\it functional variations} introduced in \cite{N_Tiba2012}, \cite{NP_Tiba2009}
allow simultaneous topological and boundary variations. This method is of fixed domain type
and avoids drawbacks like remeshing and recomputing the mass matrix, in each iteration. In fact, in \cite{MT2019}, again for Dirichlet boundary conditions and distributed cost, we have put together all these developments and obtained a complete approximation technique with the potential to solve general shape optimization problems (general cost functionals, general boundary conditions, various differential operators, including parabolic operators as well, etc.). We continue in this paper with the case of boundary observation and we show that the new approach, with certain natural modifications and adaptations, gives good results too. 
Notice that such ideas are also applicable in free boundary problems, for instance for
fluid-structure interaction \cite{HMT2016}, \cite{HMT2018}.
Other applications are in optimization and optimal control \cite{Tiba2020}.

In the next section, we collect some preliminaries on the implicit parametrization
theorem, especially in dimension $d=2$, which is a case of interest in shape optimization.
The formulation of the problem is also discussed.
The approximation via penalization and its differentiability properties are analyzed
in Section 3. Next, we investigate the discretization process in Section 4.
The last section is devoted to numerical experiments.

\section{Preliminaries and problem formulation}
\setcounter{equation}{0}

In this paper, we fix our attention on the problem ($\mathcal{P}$):
\begin{equation}\label{2.1}
\min_{\Omega\in\mathcal{O}} \int_{\partial\Omega} j\left(\mathbf{x}, \nabla y_\Omega(\mathbf{x})\right)d\sigma,
\end{equation} 
subject to (\ref{1.2})-(\ref{1.3}) and with $j:D\times\mathbb{R}^2\rightarrow \mathbb{R}$, a
Caratheodory mapping. The dependence of $j$ on $y_\Omega$ is not necessary here since 
$y_\Omega = 0$ on $\partial\Omega$. A classical example is the normal derivative
$j\left(\mathbf{x}, \nabla y_\Omega(\mathbf{x})\right)
=\left| \frac{\partial y_\Omega}{\partial \mathbf{n}}\right|^2$.

According to the functional variations approach, introduced in \cite{N_Tiba2012},
\cite{NP_Tiba2009},
we consider that the family $\mathcal{O}$ of admissible domains given in (\ref{2.1}), is defined 
starting from  a family of admissible function $\mathcal{F}\subset\mathcal{C}(\overline{D})$
(where $D$ is a bounded domain in $\mathbb{R}^2$) via the relation:
\begin{equation}\label{2.2}
\Omega=\Omega_g=int\left\{ \mathbf{x}\in D;\ g(\mathbf{x}) \leq 0\right\},\quad g\in \mathcal{F}.
\end{equation} 
While relation (\ref{2.2}) defines a family of open sets (not necessarily connected), by imposing
further natural geometric constraints, relation (\ref{2.2}) defines a family of domains.
One example is the selection of the connected component containing $E$
\begin{equation}\label{2.3}
E \subset \Omega,\quad \forall \Omega \in \mathcal{O},
\end{equation} 
where $E$ is a given subdomain such that $\overline{E}\subset D$. In the formulation (\ref{2.2}),
inclusion (\ref{2.3}) is expressed as 
\begin{equation}\label{2.4}
g(\mathbf{x})\leq 0,\quad \forall \mathbf{x}\in \overline{E}.
\end{equation}
Another example is the selection of the connected component via
$\mathbf{x}^0\in \partial\Omega$, for any $\Omega$ in $\mathcal{O}$.
This can be reformulated as
\begin{equation}\label{2.5}
g(\mathbf{x}^0)=0,\quad \forall g\in \mathcal{F}
\end{equation}
if $ \mathcal{F}\subset \mathcal{C}^1(\overline{D})$ and satisfies the following conditions
(according to \cite{Tiba2018a}):
\begin{eqnarray}
g( \mathbf{x}) & > & 0,\quad\hbox{on }\partial D,
\label{2.6}\\
|\nabla g( \mathbf{x})| & > &  0,\quad\hbox{on }
\mathcal{G}=\left\{  \mathbf{x}\in D;\ g(\mathbf{x})=0\right\}.
\label{2.7}
\end{eqnarray}
This is due to the implicit functions theorem applied to the equation $g(\mathbf{x}) =0$,
around $\mathbf{x}^0$ from (\ref{2.5}).
By (\ref{2.6}), (\ref{2.7}), we get that $\Omega_g\cap\partial D=\emptyset $ for
any $g\in\mathcal{F}$ and (\ref{2.2}) can be equivalently expressed as
\begin{equation}\label{2.8}
\Omega_g=\left\{ \mathbf{x}\in D;\ g(\mathbf{x}) < 0\right\},\quad g\in \mathcal{F}.
\end{equation}

Similarly, if we want that a given manifold $\mathcal{C}\subset D$ is contained in
$\partial\Omega_g$ for any $g\in\mathcal{F}$, then we impose
\begin{equation}\label{2.9}
g(\mathbf{x})=0,\quad \mathbf{x}\in\mathcal{C},\  g\in \mathcal{F}.
\end{equation} 

We notice that the family $\mathcal{F}$ is very large and very flexible in 
imposing various geometric constraints on the admissible domains $\mathcal{O}$, via simple conditions on $\mathcal{F}$.
It includes, for instance, multimodal functions of class $\mathcal{C}^1(\overline{D})$
that may have unbounded many extremal points in $D$.
Moreover, the obtained domains $\Omega_g$ are connected but not simply connected,
in general. Consequently, our approach, allows topological optimization and
performs, in fact, simultaneous topological and boundary variations, which is a
characteristic of functional variations \cite{N_Tiba2012}, \cite{NP_Tiba2009} .
We ask that $D\subset \mathbb{R}^2$, which is an important case in shape optimization.
This restriction is due to the use of 
Poincar\'e-Bendixson type arguments, in some of the following results (see Hirsch, Smale and Devaney \cite{Hirsch2014}
, Ch. 10 or Pontryagin \cite{Pon}).

\begin{proposition}(Tiba \cite{Tiba2018a})\label{prop:2.1}
If $D\subset \mathbb{R}^2$, $\mathcal{F}\subset \mathcal{C}^2(\overline{D})$ and assumptions 
(\ref{2.6}), (\ref{2.7}) are valid, then 
$\mathcal{G}=\left\{  \mathbf{x}\in D;\ g(\mathbf{x})=0\right\}$ 
is a finite union of disjoint closed curves of class $\mathcal{C}^2$, without self intersections,
and not intersecting $\partial D$.
They are parametrized by the solution of the Hamiltonian system:
\begin{eqnarray}
x_1^\prime(t) & = & -\frac{\partial g}{\partial x_2} \left( x_1(t), x_2(t) \right), \ t\in I,
\label{2.10}\\
x_2^\prime(t) & = &  \frac{\partial g}{\partial x_1} \left( x_1(t), x_2(t) \right), \ t\in I,
\label{2.11}\\
\left(x_1(0),x_2(0)\right)& = & \mathbf{x}^0  = \left(x_1^0, x_2^0\right) \in D,
\label{2.12}
\end{eqnarray}
where some $\mathbf{x}^0$ is chosen on each component of $\mathcal{G}$.
\end{proposition}

Here, the constraint (\ref{2.4}) is not necessarily valid and $\Omega_g$
from (\ref{2.8}) is a finite union of domains, that may be multiply connected.
The existence interval $I$ from (\ref{2.10})-(\ref{2.12}) may be taken $I=\mathbb{R}$ or
just the corresponding period (the solutions of (\ref{2.10})-(\ref{2.12}) are periodic - this is the consequence of the Poincar\'e-Bendixson result and hypotheses (\ref{2.6}), (\ref{2.7})).
In higher dimension, iterated Hamiltonian systems have to be used and their solution
may be just a local one, Tiba \cite{Tiba2018}. This is the case of the implicit
parametrization method, a recent extension of the implicit function theorem.

Consider now another mapping $h\in \mathcal{C}^2(\overline{D})$ and satisfying 
(\ref{2.6}), (\ref{2.7}). We define the functional perturbation
$g+\lambda h$, $\lambda\in \mathbb{R}$ ``small'', such that (\ref{2.6}), (\ref{2.7}) are
still satisfied by $g+\lambda h$, due to some simple argument based
on the Weierstrass theorem.

\begin{proposition}(Tiba \cite{Tiba2018a})\label{prop:2.2}
If $\epsilon>0$ is small enough, there is $\lambda ( \epsilon ) > 0$ such that,
for $\lambda\in \mathbb{R}$, $|\lambda| < \lambda(\epsilon)$,
we have that $\mathcal{G}_\lambda$ in included in $V_\epsilon$
and $\mathcal{G}_\lambda$ is a finite union of $\mathcal{C}^2$ curves.
\end{proposition}

Here
\begin{eqnarray*}
\mathcal{G}_\lambda&=&\left\{ \mathbf{x}\in D;\ (g+\lambda h)(\mathbf{x})=0\right\},\\
V_\epsilon&=&\left\{ \mathbf{x}\in D;\ d[\mathbf{x},\mathcal{G}] < \epsilon\right\}
\end{eqnarray*}
with $d[\mathbf{x},\mathcal{G}]$ being the distance between a point and $\mathcal{G}$.
In particular, Proposition \ref{prop:2.2} shows that $\mathcal{G}_\lambda\rightarrow \mathcal{G}$
in the Hausdorff-Pompeiu sense, Neittaanm\"aki et al. \cite{NS_Tiba2006}, Appendix 3.

\begin{proposition}(Murea and Tiba \cite{MT2019})\label{prop:2.3}
Denote by $T_g$, $T_{\lambda}$ the periods of the Hamiltonian system (\ref{2.10})-(\ref{2.12}), 
respectively the perturbed Hamiltonian system.
Then $T_{\lambda} \rightarrow T_g$ as $\lambda \rightarrow 0$.
\end{proposition}

\begin{remark}\label{rem:2.4}
A discussion of the dependence of the period $T_g$ with respect to certain
perturbations can be found in Teschl \cite{Teschl2010}, Ch. 12.
 In general, the perturbation of a periodic
system may not be periodic and the approximation properties have an asymptotic
character, Sideris \cite{Sideris2013}. In \cite{MT2020}, we have established that the period $T_g$ has even differentiability properties with respect to functional variations and this will be used in the next Section. 
\end{remark}

\section{Approximation and differentiability}
\setcounter{equation}{0}

We shall use a variant of the penalization method from Tiba \cite{Tiba2018a},
that has good differentiability properties as well. The main new ingredient in this approach is that we penalize directly the cost functional and not the state equation as in \cite{Tiba2018b},
\cite{Tiba2016}, \cite{MT2019a}. This appears as the application of classical optimization techniques and its advantage is the possibility to extend it to any boundary conditions. We underline that the Hamiltonian handling of the unknown geometries plays an essential role in the formulation below.

\vspace{3mm}

\noindent
The penalized optimization problem is given by
\begin{equation}\label{3.1}
\min_{g,u} 
\int_{I_g}\left[ j\left(\mathbf{z}_g(t), \nabla y(\mathbf{z}_g(t))\right)
+\frac{1}{\epsilon}\left(y(\mathbf{z}_g(t))\right)^2 \right] |\mathbf{z}_g^\prime(t)|dt
\end{equation}
\begin{eqnarray}
-\Delta y & = & f+ g_+^2 u,\quad\hbox{in } D,
\label{3.2}\\
y & = & 0,\quad\hbox{on } \partial D,
\label{3.3}\\
g(\mathbf{x}) & \leq & 0,\quad\hbox{on }\overline{E}\subset D, \hbox{ given}
\label{3.4}
\end{eqnarray}
where $\mathbf{z}_g:I_g\rightarrow D $, $\mathbf{z}_g\in\left(\mathcal{C}^1(I_g)\right)^2$
is the solution of the Hamiltonian system (\ref{2.10})-(\ref{2.12}) associated to
$g\in \mathcal{F}$ and $I_g=[0,T_g]$ is its period.
In case $\partial \Omega_g$ has several components (their number is finite according
to Section 2), then the penalization part in the functional (\ref{3.1}) has to be understood as a finite sum
of terms corresponding to each component. Notice that the corresponding periods
and the initial conditions (\ref{2.12}) can be obtained via standard numerical
methods in the examples, see Remark \ref{rem:4.1}.

The minimization is performed over $g\in \mathcal{F}$, satisfying (\ref{3.4}), (\ref{2.6}),
(\ref{2.7}) and $u$ measurable such that $g_+^2 u\in L^p(D)$, $p>2$.
It is possible that the original cost (\ref{2.1}) (the first term in (\ref{3.1})) is
defined just on one component of $\partial\Omega_g$ and this can be singled out by a condition
like (\ref{2.5}) and a corresponding given $\mathbf{x}^0\notin \overline{E}$.
However the penalization term in (\ref{3.1}) has to be defined on all the
components of $\partial\Omega_g$ since it controls in fact the Dirichlet condition
(\ref{1.3}). For simplicity, we shall not investigate such details here, related to (\ref{3.1}).

If $\partial D$ is in $\mathcal{C}^{1,1}$, then the state $y\in W^{2,p}(D)\cap H^1_0(D)$,
due to (\ref{3.2}), (\ref{3.3}). Consequently $y\in \mathcal{C}^1(\overline{D})$.
Then, the cost functionals (\ref{2.1}), (\ref{3.1}) make sense since $\nabla y$
is continuous in $\overline{D}$ and similar regularity properties are valid on $\Omega_g$ under the assumptions on $g\in \mathcal{F}$.

\begin{proposition}\label{prop:3.1}
Let $j(\cdot,\cdot)$ be a Carath\'eodory function on $D\times\mathbb{R}^2$,
bounded by a constant from below. Let $[y_n^\epsilon, g_n^\epsilon, u_n^\epsilon]$ be
a minimizing sequence in
the penalized problem (\ref{3.1})-(\ref{3.4}), for some given $\epsilon>0$. 
Then, on a subsequence denoted by $n(m)$ the (not necessarily admissible) pairs
$[\Omega_{g_{n(m)}^\epsilon}, y_{n(m)}^\epsilon]$ give a minimizing cost in (\ref{2.1}),
satisfy (\ref{1.2}) in $\Omega_{g_{n(m)}^\epsilon}$ and (\ref{1.3}) is fulfilled with
a perturbation of order $\epsilon^{1/2}$ on $\partial \Omega_{g_{n(m)}^\epsilon}$.
\end{proposition}

\noindent
\textbf{Proof.} 
Let $[y_{g_m}, g_m]\in W^{2,p}(\Omega_{g_{m}})\times \mathcal{F}$ be a minimizing sequence in the problem
(\ref{2.1}), (\ref{1.2}), (\ref{1.3}), (\ref{3.4}) where $\Omega=\Omega_g$ is defined by (\ref{2.8})
and $g$ satisfies $g_+^2 u\in L^p(D)$.
By Proposition \ref{prop:2.1}, $\partial \Omega_g$ is of class $\mathcal{C}^2$ and this ensures 
the regularity for (\ref{1.2}), (\ref{1.3}) since $f\in L^p(D)$.

Take $\widetilde{y}_{g_{m}}\in W^{2,p}(D\setminus \overline{\Omega_{g_{m}}})$, not unique, given by the trace theorem
such that $\widetilde{y}_{g_{m}}=y_{g_m}$ on $\partial \Omega_{g_{m}}$, 
$\frac{\partial \widetilde{y}_{g_{m}}}{\partial \mathbf{n}} = \frac{\partial y_{g_{m}}}{\partial \mathbf{n}}$
on $\partial \Omega_{g_{m}}$, $\widetilde{y}_{g_{m}}=0$ on $\partial D$.
We define an admissible control $u_{g_{m}}$ in (\ref{3.2}) by 
\begin{equation}
u_{g_{m}} = -\frac{\Delta \widetilde{y}_{g_{m}} + f}{(g_m)_+^2},
\quad\hbox{in }  D\setminus \overline{\Omega}_{g_{m}}
\end{equation}
and zero otherwise.
It yields $(g_m)_+^2u_{g_{m}}\in L^p(D)$ and this control pair is admissible for the
problem (\ref{3.1})-(\ref{3.4}). Moreover, the corresponding state 
$\overline{y}_{g_{m}}$ in (\ref{3.2})-(\ref{3.3}) is obtained by concatenation of $y_{g_{m}}$ and $\widetilde{y}_{g_{m}}$ and
the associated penalization term in (\ref{3.1}) is null, due to (\ref{1.3}).

We get the inequality:
\begin{eqnarray}
&&\int_{I_{g_{n(m)}^\epsilon}} 
\left[ j\left( \mathbf{z}_{g_{n(m)}^\epsilon}, \nabla y_{n(m)}^\epsilon(\mathbf{z}_{g_{n(m)}^\epsilon})\right)
+\frac{1}{\epsilon}\left(y_{n(m)}^\epsilon
(\mathbf{z}_{g_{n(m)}^\epsilon})\right)^2\right]
|\mathbf{z}_{g_{n(m)}^\epsilon}^\prime|dt
\nonumber
\\
&\leq &
\int_{\partial \Omega_{g_m}} j\left(\mathbf{x}, \nabla y_{g_{m}}(\mathbf{x})\right)d\sigma
\rightarrow \inf (\mathcal{P}),
\label{3.6}
\end{eqnarray}
for $n(m)$ big enough, due to the minimizing property of the sequence
$[y_n^\epsilon, g_n^\epsilon, u_n^\epsilon]$, respectively $[y_{g_m}, g_m]$. 
By (\ref{3.6}) we infer
\begin{equation}\label{3.7}
\int_{\partial \Omega_{g_{m}^\epsilon}} \left( y_{n(m)}^\epsilon \right)^2 d\sigma
\leq C \epsilon
\end{equation}
with $C$ a constant independent of $ \epsilon $,
$m$ since $j$ is bounded below by a constant. 
Relation (\ref{3.7}) proves the last statement in the proposition.
As $(g_{n(m)}^\epsilon)_+$ is null in $\Omega_{g_{n(m)}^\epsilon}$, we see that (\ref{1.2}) is
satisfied here, due to (\ref{3.2}). The minimizing property with respect to the original
cost (\ref{2.1}) is a clear consequence of (\ref{3.6}).\quad$\Box$

\begin{remark}\label{rem:3.2}
In  \cite{Tiba2018a}, \cite{MT2020} a detailed study of the approximating properties
with respect to $\epsilon\rightarrow 0$, is performed in related problems.
\end{remark}

We consider now $[u,g]\in L^p(D)\times\mathcal{F}$, $p>2$, satisfying (\ref{3.4}),
(\ref{2.5}) together with perturbations $[u+\lambda v,g+\lambda r]$,
$\lambda \in \mathbb{R}$, $v\in L^p(D)$, such that (\ref{3.4}),
(\ref{2.5}) are satisfied by $r\in \mathcal{F}$. The  state system is,
in fact, given by (\ref{3.2}), (\ref{3.3}), (\ref{2.10})-(\ref{2.12}) and the corresponding perturbed system has solutions $y^\lambda, \mathbf{z}_{g+\lambda r}$. We study
its differentiability properties.

\begin{proposition}\label{prop:3.3}
The system in variations corresponding to (\ref{3.2}), (\ref{3.3}),
(\ref{2.10})-(\ref{2.12}) is:
\begin{eqnarray}
-\Delta q & = & g_+^2 v + 2g_+ u\,r,\quad\hbox{in } D,
\label{3.8}\\
q & = & 0,\quad\hbox{on } \partial D,
\label{3.9}\\
w_1^\prime & = & -\nabla\partial_2 g(\mathbf{z}_g) \cdot \mathbf{w}
- \partial_2 r(\mathbf{z}_g), \hbox{ in }I_g,
\label{3.10}\\
w_2^\prime & = & \nabla\partial_1 g(\mathbf{z}_g) \cdot \mathbf{w}
+ \partial_1 r(\mathbf{z}_g), \hbox{ in }I_g,
\label{3.11}\\
w_1(0)& = &0,\ w_2(0) = 0,
\label{3.12}
\end{eqnarray}
where $q=\lim_{\lambda\rightarrow 0}\frac{y^\lambda-y}{\lambda}$, 
$\mathbf{w}=[w_1,w_2]=\lim_{\lambda\rightarrow 0}
\frac{\mathbf{z}_{g+\lambda r} - \mathbf{z}_g}{\lambda}$
with $y^\lambda\in W^{2,p}(D)\cap H^1_0(D)$ being the solution of (\ref{3.2}), (\ref{3.3})
corresponding to $g+\lambda r$, $u+\lambda v$ and
``$\cdot$'' is the scalar product in $\mathbb{R}^2$.
The limits exist in the spaces of $y$, $\mathbf{z}_g$, respectively. 
\end{proposition}

\noindent
\textbf{Proof.} 
Subtracting the equations of $y^\lambda$ (i.e. (\ref{3.2}), (\ref{3.3}) with perturbed
controls) and $y$, we get 
\begin{equation}\label{3.13}
-\Delta  \frac{y^\lambda-y}{\lambda}
= \frac{1}{\lambda}\left[ (g+\lambda r)_+^2(u+\lambda v)
- g_+^2 u\right],\quad\hbox{in } D,
\end{equation}
with zero boundary conditions on $\partial D$.
A standard passage to the limit in (\ref{3.13}), gives (\ref{3.8}), (\ref{3.9}).

For (\ref{3.10})-(\ref{3.12}), the argument is similar as in Proposition 6,
Tiba \cite{Tiba2013}.  The convergence is in $\mathcal{C}^1(I_g)$ on the whole
sequence $\lambda\rightarrow 0$ due to the uniqueness property for the linear systems
(\ref{3.8})-(\ref{3.12}) and the periodicity of the solutions $\mathbf{z}_g$,
$\mathbf{z}_{g+\lambda r}$ by Proposition \ref{prop:2.1}.\quad$\Box$

We assume now that $j(\mathbf{x},\cdot)$ is $\mathcal{C}^1 (\mathbb{R}^2)$, $j(\mathbf{x}^0,\cdot)\equiv 0$
and $f\in W^{1,p}(D)$, $\partial D$ is in $\mathcal{C}^{2,1}$. Notice that by imposing
$\mathcal{F}\subset \mathcal{C}^2(\overline{D})$, we get that $g_+^2\in W^{1,\infty}(D)$
and $g_+^2u\in W^{1,p}(D)$ if $u\in W^{1,p}(D)$.

\begin{proposition}\label{prop:3.4}
Under the above hypotheses, if $y(\mathbf{x}^0)=0$, then
the directional derivative of the penalized cost (\ref{3.1}) in the direction
$[v,r]\in W^{1,p}(D)\times \mathcal{F}$ is given by:
\begin{eqnarray}
&&
\int_{I_g} \nabla_1 j\left(\mathbf{z}_g(t), \nabla y(\mathbf{z}_g(t))\right)  
\cdot \mathbf{w}(t)\, |\mathbf{z}_g^\prime(t)| dt
\nonumber\\
& + & 
\int_{I_g} \nabla_2 j\left(\mathbf{z}_g(t), \nabla y(\mathbf{z}_g(t))\right)  
\cdot H\left(y(\mathbf{z}_g(t))\right) \cdot \mathbf{w}(t)\, |\mathbf{z}_g^\prime(t)| dt
\nonumber\\
& + & 
\int_{I_g} \nabla_2 j\left(\mathbf{z}_g(t), \nabla y(\mathbf{z}_g(t))\right)  
\cdot \nabla q(\mathbf{z}_g(t))\, |\mathbf{z}_g^\prime(t)| dt
\nonumber\\
& + & \frac{2}{\epsilon}\int_{I_g} y(\mathbf{z}_g(t))
\left[\nabla y(\mathbf{z}_g(t))\cdot \mathbf{w}(t)
+q(\mathbf{z}_g(t)) \right]\, |\mathbf{z}_g^\prime(t)| dt
\nonumber\\
&+&\int_{I_g}
\left[ j\left(\mathbf{z}_g(t), \nabla y(\mathbf{z}_g(t))\right) 
+\frac{1}{\epsilon}\left(y(\mathbf{z}_g(t))\right)^2\right]
\frac{\mathbf{z}_g^\prime(t)\cdot \mathbf{w}^\prime(t)}{|\mathbf{z}_g^\prime(t)|} dt .
\label{3.14}
\end{eqnarray}

\end{proposition}

\noindent The notations are explained in the proof.

\noindent
\textbf{Proof.} 
We compute
\begin{eqnarray}
&&\lim_{\lambda\rightarrow 0}\frac{1}{\lambda}\left\{ \int_{I_{g+\lambda r}}
\left[ j\left(\mathbf{z}_{g+\lambda r}(t), \nabla y^\lambda(\mathbf{z}_{g+\lambda r}(t))\right) 
+\frac{1}{\epsilon}\left(y^\lambda(\mathbf{z}_{g+\lambda r}(t))\right)^2\right]
\, |\mathbf{z}_{g+\lambda r}^\prime(t)| dt \right.
\nonumber\\
&-&\left. \int_{I_g}
\left[ j\left(\mathbf{z}_g(t), \nabla y(\mathbf{z}_g(t))\right) 
+\frac{1}{\epsilon}\left(y(\mathbf{z}_g(t))\right)^2\right]\, |\mathbf{z}_g^\prime(t)| dt\right\},
\label{3.15}
\end{eqnarray}
where we use the notations from Proposition \ref{prop:3.3}. The above assumptions
on $\mathcal{F}$, $u$, $v$ ensure that $y^\lambda, y\in W^{3,p}(D)\subset \mathcal{C}^2(\overline{D})$, 
Grisvard \cite{Grisvard1985}, and $y^\lambda\rightarrow y$ in $\mathcal{C}^2(\overline{D})$,
$\mathbf{z}_{g+\lambda r}\rightarrow\mathbf{z}_g$ in $\mathcal{C}^2(I_g)$.

We study first the term:
\begin{eqnarray}
&&\frac{1}{\lambda}\int^{T_{g+\lambda r}}_{T_g}
\left[
j\left(\mathbf{z}_{g+\lambda r}(t), \nabla y^\lambda(\mathbf{z}_{g+\lambda r}(t))\right) 
+\frac{1}{\epsilon}\left(y^\lambda(\mathbf{z}_{g+\lambda r}(t))\right)^2\right]
\, |\mathbf{z}_{g+\lambda r}^\prime(t)| dt 
\nonumber\\
&=&\frac{T_{g+\lambda r} - T_g}{\lambda}
\left[
j\left(\mathbf{z}_{g+\lambda r}(\tau), \nabla y^\lambda(\mathbf{z}_{g+\lambda r}(\tau))\right) 
+\frac{1}{\epsilon}\left(y^\lambda(\mathbf{z}_{g+\lambda r}(\tau))\right)^2\right]
\, |\mathbf{z}_{g+\lambda r}^\prime(\tau)|\nonumber\\
&\rightarrow& 0
\label{3.16}
\end{eqnarray}
due to the differentiability properties of $T_g$ with respect to functional variations $g+\lambda r$ (see  \cite{MT2020}) and the convergence properties of $y^\lambda$,
$\mathbf{z}_{g+\lambda r}$ and the regularity assumptions on $j(\cdot,\cdot)$.
In  (\ref{3.16})  $\tau$, is some intermediary 
point in $[T_g, T_{g+\lambda r}]$ depending on $g$, $r$, $\lambda$.
The assumptions on $j(\mathbf{x}^0,\cdot)$ and $y(\mathbf{x}^0)$ give that the term studied in 
(\ref{3.16}) has null limit and can be neglected.

For the first term in (\ref{3.15}), we have
\begin{eqnarray*}
&& \int_{I_g}\frac{j\left(\mathbf{z}_{g+\lambda r}(t), \nabla y^\lambda(\mathbf{z}_{g+\lambda r}(t))\right)
-j\left(\mathbf{z}_g(t), \nabla y(\mathbf{z}_g(t))\right) }{\lambda}
|\mathbf{z}_{g}^\prime(t)|dt\\
&\rightarrow&
\int_{I_g}\nabla_1 j\left(\mathbf{z}_g(t), \nabla y(\mathbf{z}_g(t))\right)  
\cdot \mathbf{w}(t)\, |\mathbf{z}_g^\prime(t)| dt\\
&+&\int_{I_g} \nabla_2 j\left(\mathbf{z}_g(t), \nabla y(\mathbf{z}_g(t))\right)  
\cdot H\left(y(\mathbf{z}_g(t))\right) \cdot \mathbf{w}(t)\, |\mathbf{z}_g^\prime(t)| dt\\
&+&\int_{I_g} \nabla_2 j\left(\mathbf{z}_g(t), \nabla y(\mathbf{z}_g(t))\right)  
\cdot \nabla q(\mathbf{z}_g(t))\, |\mathbf{z}_g^\prime(t)| dt
\end{eqnarray*}
where $\nabla_1$, $\nabla_2$ denote the gradient of $j()$ with respect to the first two arguments,
respectively the last two arguments, $H(y)$ is the Hessian matrix and $\mathbf{w}$ is given
by (\ref{3.10})-(\ref{3.12}), $q$ is given by (\ref{3.8}), (\ref{3.9}).

Consider now a second part from (\ref{3.15}):
\begin{eqnarray*}
&& \lim_{\lambda\rightarrow 0}\frac{1}{\lambda\epsilon} \int_{I_g}
\left[ \left(y^\lambda(\mathbf{z}_{g+\lambda r}(t))\right)^2
-\left(y(\mathbf{z}_g(t))\right)^2\right]\, |\mathbf{z}_g^\prime(t)| dt\\
&=&\frac{2}{\epsilon}\int_{I_g}
y(\mathbf{z}_g(t))\left[ 
\nabla y(\mathbf{z}_g(t))\cdot\mathbf{w}(t) + q(\mathbf{z}_g(t))
\right]
\, |\mathbf{z}_g^\prime(t)| dt.
\end{eqnarray*}

It remains to complete:
\begin{eqnarray*}
&& \lim_{\lambda\rightarrow 0}\frac{1}{\lambda} \int_{I_g}
\left[
j\left(\mathbf{z}_{g+\lambda r}(t), \nabla y^\lambda(\mathbf{z}_{g+\lambda r}(t))\right)
+\frac{1}{\epsilon}\left(y^\lambda(\mathbf{z}_{g+\lambda r}(t))\right)^2
\right]\left( |\mathbf{z}_{g+\lambda r}^\prime(t)| - |\mathbf{z}_g^\prime(t)|\right) dt
\\
&&=\int_{I_g}
\left[ j\left(\mathbf{z}_g(t), \nabla y(\mathbf{z}_g(t))\right) 
+\frac{1}{\epsilon}\left(y(\mathbf{z}_g(t))\right)^2\right]
\frac{\mathbf{z}_g^\prime(t)\cdot \mathbf{w}^\prime(t)}{|\mathbf{z}_g^\prime(t)|} dt.
\end{eqnarray*}

Notice that $|\mathbf{z}_g^\prime(t)|\neq 0$ on $I_g$ due to  (\ref{2.10})-(\ref{2.12})
and (\ref{2.7}). The above computations are based on appropriate interpolation of terms
and differentiability properties of the involved quantities. In particular, 
in the last computation, the critical case is avoided.\quad$\Box$

Denote by
$A:\mathcal{C}^2(\overline{D})\times  W^{1,p}(D)\rightarrow W^{3,p}(D)\cap W_0^{1,p}(D)$
the linear continuous operator 
$r, v\rightarrow q$ given by (\ref{3.8}), (\ref{3.9}) and by
$B: \mathcal{C}^2 \left( \overline{D} \right) \rightarrow \mathcal{C}^1(I_g)^2 $
the linear continuous operator given by (\ref{3.10})-(\ref{3.12}),
via the relation $r\rightarrow \mathbf{w}$. In these definitions, $g\in \mathcal{C}^2(\overline{D})$ and
$u\in W^{1,p}(D)$ are fixed.

\begin{corollary}\label{cor:3.6}
The relation (\ref{3.14}) can be rewritten as:
\begin{eqnarray}
&&
\int_{I_g} \nabla_1 j\left(\mathbf{z}_g(t), \nabla y(\mathbf{z}_g(t))\right)  
\cdot Br(\mathbf{z}_g(t))\, |\mathbf{z}_g^\prime(t)| dt
\nonumber\\
& + & 
\int_{I_g} \nabla_2 j\left(\mathbf{z}_g(t), \nabla y(\mathbf{z}_g(t))\right)  
\cdot H\left(y(\mathbf{z}_g(t))\right) \cdot  Br(\mathbf{z}_g(t))\, |\mathbf{z}_g^\prime(t)| dt
\nonumber\\
& + & 
\int_{I_g} \nabla_2 j\left(\mathbf{z}_g(t), \nabla y(\mathbf{z}_g(t))\right)  
\cdot \nabla A(r,v)(\mathbf{z}_g(t))\, |\mathbf{z}_g^\prime(t)| dt
\nonumber\\
& + & \frac{2}{\epsilon}\int_{I_g} y(\mathbf{z}_g(t))
\left[\nabla y(\mathbf{z}_g(t)) \cdot Br(\mathbf{z}_g(t))
+ A(r,v)(\mathbf{z}_g(t)) \right]\, |\mathbf{z}_g^\prime(t)| dt
\nonumber\\
&+&\int_{I_g}
\left[ j\left(\mathbf{z}_g(t), \nabla y(\mathbf{z}_g(t))\right) 
+\frac{1}{\epsilon}\left(y(\mathbf{z}_g(t))\right)^2\right]
\frac{\mathbf{z}_g^\prime(t)}{|\mathbf{z}_g^\prime(t)|} 
\cdot [-\partial_2 r, \partial_1 r](\mathbf{z}_g(t)) dt
\nonumber\\
&+&
\int_{I_g}C(t)\cdot Br(\mathbf{z}_g(t)) dt.
\label{3.17}
\end{eqnarray}
\end{corollary}

Here $C(t)$ is a vector obtained by replacing $\mathbf{w}^\prime(t)$ as expressed in
(\ref{3.10}), (\ref{3.11}) and separating the part including $[-\partial_2 r, \partial_1 r]$.

\begin{remark}\label{rem:3.5}
The regularity hypotheses are natural and necessary when making variations of boundary
integrals. The conditions $j(\mathbf{x}^0,\cdot)\equiv 0$ can be obtained by a translation and
$y(\mathbf{x}^0)=0$  reflects that $\mathbf{x}^0 \in \mathcal{G}$ and the admissible states in the original
shape optimization problem are automatically null on $\mathcal{G}$. It is possible to remove these two conditions (see  \cite{MT2020}), but the relation (\ref{3.14}) becomes more complex.
\end{remark}

\section{Finite element discretization}
\setcounter{equation}{0}

We assume that $D$ is polygonal and let $\mathcal{T}_h$ be a triangulation of $D$
where $h$ is the size of $\mathcal{T}_h$.
We introduce the linear space 
$$
\mathbb{W}_h = \left\{ \varphi_h \in \mathcal{C}(\overline{D}) ;
\ {\varphi_h}_{|T} \in \mathbb{P}_3(T),\ \forall T \in \mathcal{T}_h  \right\}
$$
where $\mathbb{P}_3$ is the piecewise cubic finite element.
We use a standard basis of $\mathbb{W}_h$, $\left\{ \phi_i\right\}_{i\in I}$, where
$I=\{1, \dots, n\}$ and
$\phi_i$ is the hat function
associated to the node $A_i$, see for example \cite{Ciarlet2002}, \cite{Raviart2004}.
There are ten nodes for the cubic finite element on a triangle.

We can approach $g$ and $u$ by the finite element functions
$g_h=\sum_{i\in I} G_i \phi_i$ and 
$u_h=\sum_{i\in I} U_i \phi_i$. We introduce the $\mathbb{R}^n$ vectors
$G=(G_i)_{i\in I}^T$, $U=(U_i)_{i\in I}^T$ and $g_h$ can be identified by $G$, etc.
It is possible to use for $u$ a low order finite element,
like piecewise linear $\mathbb{P}_1$.

We also set
$$
\mathbb{V}_h = \left\{ \varphi_h \in \mathbb{W}_h ;
\ \varphi_h = 0\hbox{ on }\partial D \right\},
$$
$I_0=\{i\in I;\ A_i\notin \partial D \}$ where $n_0=card(I_0)$ and the vector
$$
F=\left( F_i \right)_{i\in I_0} = \left(\int_D f\phi_i d\mathbf{x}\right)_{i\in I_0}\in\mathbb{R}^{n_0}.
$$
The discrete weak formulation of (\ref{3.2})-(\ref{3.3}) is: 
find $y_h \in \mathbb{V}_h $ such that
\begin{equation}
\label{4.1}
\int_D \nabla y_h \cdot\nabla  \varphi_h\, d\mathbf{x}
= \int_D \left(f_h + (g_h)_+^2u_h\right) \varphi_h\, d\mathbf{x},
\quad \forall \varphi_h \in \mathbb{V}_h.
\end{equation}
The finite element approximations of $y$ is $y_h(\mathbf{x})=\sum_{j\in I_0} Y_j \phi_j(\mathbf{x})$ 
with $Y=(Y_j)_{j\in I_0}^T\in\mathbb{R}^{n_0}$ and similarly for $f$, $f_h$, $F$.

Let us define $K$ the square matrix of order $n_0$ by
$$
K=(K_{ij})_{i\in I_0,j\in I_0}, \quad
K_{ij} =  \int_D \nabla \phi_j \cdot \nabla \phi_i d\,\mathbf{x}
$$
and the $n_0\times n$ matrix $B^1(G)$ defined by
$$
B^1(G)=(B_{ij}^1)_{i\in I_0,j\in I},\quad 
B_{ij}^1=  \int_D  (g_h)_+^2 \phi_j\phi_i d\mathbf{x}.
$$
The matrix $K$ is symmetric, positive definite and
the linear system associated to the state system (\ref{3.2})-(\ref{3.3}) is:
\begin{equation}\label{4.2}
KY=F+B^1(G)U.
\end{equation}

For the time step $\Delta t >0$, the forward Euler scheme can be used:
\begin{eqnarray}
Z_{k+1}^1 & = & Z_{k}^1 - \Delta t \frac{\partial g_h}{\partial x_2}\left(Z_{k}^1, Z_{k}^2\right),
\label{4.3}\\
Z_{k+1}^2 & = & Z_{k}^2 + \Delta t \frac{\partial g_h}{\partial x_1}
\left( Z_{k}^1, Z_{k}^2 \right) ,
\label{4.4}\\
(Z_0^1, Z_0^2) & = & \left( x_1^0, x_2^0 \right) ,
\label{4.5}
\end{eqnarray}
for $k=0,1,\dots ,$ in order to solve numerically the ODE system (\ref{2.10})-(\ref{2.12}).
We set $Z_k=(Z_k^1,Z_k^2)$, in fact, $Z_k$ is an approximation of $\mathbf{z}_g(t_k)$,
where $t_k=k\Delta t$, $k\in\mathbb{N}$.
When $Z_m$, for some $m\in \mathbb{N}^*$ is ``close'' to $Z_0$,
we stop the algorithm and we set the computed period
$T_g=t_m$.
We have the uniform partition $[t_0,\dots,t_k,\dots,t_m]$ of $[0,T_g]$. We denote
$Z=(Z^1,Z^2)$ in $\mathbb{R}^m \times \mathbb{R}^m$, with
$Z^1 = (Z_k^1)_{1\leq k\leq m}^T$ and $Z^2 = (Z_k^2)_{1\leq k\leq m}^T$.
One can apply more efficient numerical methods, like explicit Runge-Kutta,
however we use (\ref{4.3})-(\ref{4.5}) for the sake of simplicity.

We define the function $Z:[0,T_g]\rightarrow \mathbb{R}^2$
$$
Z(t)=\frac{t_{k+1}-t}{\Delta t} Z_k + \frac{t-t_k}{\Delta t} Z_{k+1},\quad t_k < t \leq t_{k+1}
$$
for $k=0,1,\dots ,m-1$. We remark that $Z$ is derivable on each interval $(t_k, t_{k+1})$
and $Z^\prime(t)=\frac{1}{\Delta t}(Z_{k+1}^1-Z_k^1,Z_{k+1}^2-Z_k^2)$ for $t_k < t \leq t_{k+1}$.
We define the $n_0\times n_0$ matrix $N(Z)$ as follow
$$
N(Z) = \left( 
\int_0^{T_g} \phi_j(Z(t))\phi_i(Z(t)) |Z^\prime(t)| \, dt
\right)_{i\in I_0, j\in I_0}
$$
and, with this notation, the second term of (\ref{3.1}) is approached by
$\frac{1}{\epsilon} Y^T N(Z) Y$.

We define the partial derivatives for a piecewise cubic function.
If $g_h\in \mathbb{W}_h$ and $G\in \mathbb{R}^n$ such that
$g_h(\mathbf{x})=\sum_{i\in I} G_i \phi_i(\mathbf{x})$,
we set $\Pi_h^1G \in \mathbb{R}^n$ 
$$
\left(\Pi_h^1G\right)_i
=\frac{1}{\sum_{j\in J_i} area(T_j)}
\sum_{j\in J_i} area(T_j)\partial_1 {g_h}_{|T_j}(A_i)
$$
here $J_i$ represents the set of index $j$ such that the node $A_i$ belongs to
the triangle $T_j$.
In each triangle $T_j$, the finite element function $g_h$ is a cubic polynomial function,
then $\partial_1 {g_h}_{|T_j}$ is well defined.
In the same way, we construct $\Pi_h^2 G \in \mathbb{R}^n$ for $\partial_2$.
We have that, $\Pi_h^1$ and $\Pi_h^2$ are two square matrices of order $n$ depending
on $\mathcal{T}_h$.

We define
$$
\partial_1^h g_h(\mathbf{x})=\sum_{i\in I} \left(\Pi_h^1G\right)_i \phi_i(\mathbf{x}) \in \mathbb{W}_h
$$
and similarly for $\partial_2^h g_h$.
Putting $\nabla^h g_h = (\partial_1^h g_h, \partial_2^h g_h)$ and
since $y_h\in \mathbb{V}_h\subset \mathbb{W}_h$, we can also define $\partial_1^h y_h$
and $\partial_2^h y_h$.

A typical objective function $j$ depends on the normal derivative
$\frac{\partial y}{\partial \mathbf{n}}$.
Here, the outward unit normal vector $\mathbf{n}$ of the domain $\Omega_g$ is
approached by
\begin{equation}\label{4.6}
\mathbf{n}^h(\mathbf{x})=
\frac{1}{\sqrt{(\partial_1^h g_h(\mathbf{x}))^2 + (\partial_2^h g_h(\mathbf{x}))^2 }}
\nabla^h g_h(\mathbf{x}).
\end{equation}
The first term of (\ref{3.1}) can be approached by
$$
J_1(G,Z,Y)=\int_0^{T_g} j\left( Z(t),\nabla^h y_h(Z(t))\right) |Z^\prime(t)| dt
$$
and the discrete form of the optimization problem (\ref{3.1})-(\ref{3.3})
is
\begin{equation}\label{4.7}
\min_{G,U\in \mathbb{R}^n} J(G,U)=J_1(G,Z,Y)+\frac{1}{\epsilon}Y^TN(Z)Y
\end{equation}
subject to (\ref{4.2}). We remark that, $Y$ depends on $G$ and $U$ from (\ref{4.2})
and $Z$ depends on $G$ from (\ref{4.3})-(\ref{4.5}).
For (\ref{3.4}), we have to impose similar sign conditions on $G$.


Let $r_h$, $v_h$ be in $\mathbb{W}_h$ and $R$, $V$ in $\mathbb{R}^n$
be the associated vectors.
The discrete weak formulation of (\ref{3.8})-(\ref{3.9}) is: 
find $q_h\in \mathbb{V}_h$ such that
\begin{equation}\label{4.8}
  \int_D \nabla q_h\cdot\nabla \varphi_h \, d\mathbf{x}
  = \int_D \left( (g_h)_+^2 v_h + 2(g_h)_+ u_h r_h\right)
  \varphi_h \, d\mathbf{x},
  \quad \forall \varphi_h \in \mathbb{V}_h.
\end{equation}
We set $Q\in \mathbb{R}^{n_0}$ the vector associated to $q_h$ and we construct
the $n_0\times n$ matrix $C^1(G,U)$ defined by
$$
C^1(G,U) = \left(\int_D  2 (g_h)_+ u_h \phi_j \phi_i\, d\mathbf{x}
\right)_{i\in I_0, j\in I}.
$$
The linear system of (\ref{4.8}) is
\begin{equation}\label{4.9}
KQ=B^1(G)V + C^1(G,U)R.
\end{equation}

The term containing $q$ at the fourth line of (\ref{3.14}) is approched by
\begin{equation}
\label{4.10}
\frac{2}{\epsilon} Y^T N(Z) Q
\end{equation}
where the matrix $N(Z)$ was defined in the previous subsection.

The numerical integration over the interval $I_g$ is obtained using the right Riemann sum
\cite{wiki}.
We set $F^3=(F^3_i)\in \mathbb{R}^{n_0}$ by
$$
F^3_i=
\sum_{k=1}^m
\Delta t \nabla_2 j\left(\mathbf{Z}(t_k), \nabla^h y_h(\mathbf{Z}(t_k))\right)  
\cdot \nabla^h \phi_i(\mathbf{Z}(t_k))\, |\mathbf{Z}^\prime(t_k)|
$$
for $i\in I_0$. 
The third line of (\ref{3.14}) is approched by
\begin{equation}\label{4.11}
(F^3)^T\ Q.
\end{equation}

In order to solve the ODE system (\ref{3.10})-(\ref{3.12}), we use
the backward Euler scheme on the partition constructed before:
\begin{eqnarray}
W_{k+1}^1  & = & W_{k}^1
  - \Delta t \nabla_h \partial_2^h g_h(Z_{k+1})\cdot (W_{k+1}^1,W_{k+1}^2)\label{4.12}\\
&&
- \Delta t \partial_2^h r_h(Z_{k+1}),
\nonumber\\
 W_{k+1}^2  &= & W_{k}^2
+ \Delta t \nabla_h \partial_1^h g_h \left(Z_{k+1}\right)\cdot (W_{k+1}^1,W_{k+1}^2)\label{4.13}\\
&& + \Delta t \partial_1^h r_h\left(Z_{k+1}\right),
\nonumber\\
 W_0^1 & = & 0,\ W_0^2=0,
\label{4.14}
\end{eqnarray}
for $k=0,\dots , m-1$.
Contrary to the system (\ref{2.10})-(\ref{2.12}), the system (\ref{3.10})-(\ref{3.12}) 
is linear in $\mathbf{w}$ and we can use without difficulties an implicit method to solve it.

We set $W_k=(W_k^1,W_k^2)$ and  $W_k$ is an approximation of $\mathbf{w}(t_k)$. We write
$W=(W^1,W^2)$ in $\mathbb{R}^m \times \mathbb{R}^m$, with
$W^1 = (W_k^1)_{1\leq k\leq m}^T$ and $W^2 = (W_k^2)_{1\leq k\leq m}^T$.
The function $W:[0,T_g]\rightarrow \mathbb{R}^2$ can be constructed in the same
way as for $Z$
$$
W(t)=\frac{t_{k+1}-t}{\Delta t} W_k + \frac{t-t_k}{\Delta t} W_{k+1},\quad t_k < t \leq t_{k+1}
$$
for $k = 0, 1, \dots , m-1$. We have $W(t_k) = W_k$ and
$W^\prime(t)=\frac{1}{\Delta t}(W_{k+1}^1-W_k^1,W_{k+1}^2-W_k^2)$ for $t_k < t \leq t_{k+1}$.

We denote
\begin{eqnarray*}
\Lambda_1(t)&=& \nabla_1 j\left(\mathbf{Z}(t), \nabla^h y_h(\mathbf{Z}(t))\right)  
|\mathbf{Z}^\prime(t)| \in \mathbb{R}^2\\
\Lambda_2(t)&=&\nabla_2 j\left(\mathbf{Z}(t), \nabla^h y_n(\mathbf{Z}(t))\right)  
\cdot H^h\left(y_h(\mathbf{Z}(t))\right)|\mathbf{Z}^\prime(t)| \in \mathbb{R}^2 \\
\Lambda_4(t)&=&y_h(\mathbf{Z}(t))
\nabla^h y_h(\mathbf{Z}(t))|\mathbf{Z}^\prime(t)| \in \mathbb{R}^2
\end{eqnarray*}
and we introduce the vectors:\\
$\widetilde{\Lambda}_1=\left(\widetilde{\Lambda}_1^1,\widetilde{\Lambda}_1^2\right)
\in \mathbb{R}^m \times \mathbb{R}^m$ with the components 
$(\Delta t) \Lambda_1 (t_{k})$, $1\leq k\leq m$,\\
$\widetilde{\Lambda}_2=\left(\widetilde{\Lambda}_2^1,\widetilde{\Lambda}_2^2\right)
\in \mathbb{R}^m \times \mathbb{R}^m$ with the components 
$(\Delta t) \Lambda_2 (t_{k})$, $1\leq k\leq m$ and\\
$\widetilde{\Lambda}_4=\left(\widetilde{\Lambda}_4^1,\widetilde{\Lambda}_4^2\right)
\in \mathbb{R}^m \times \mathbb{R}^m$ with the components 
$(\Delta t) \Lambda_4 (t_{k})$, $1\leq k\leq m$.

The first, second and the term containing $\mathbf{w}$ at the fourth line of (\ref{3.14})
are approched by
\begin{equation}\label{4.15}
(\widetilde{\Lambda}_1^1)^T W^1 + (\widetilde{\Lambda}_1^2)^T W^2  
+(\widetilde{\Lambda}_2^1)^T W^1 + (\widetilde{\Lambda}_2^2)^T W^2
+\frac{2}{\epsilon}\left( (\widetilde{\Lambda}_4^1)^T W^1 + (\widetilde{\Lambda}_4^2)^T W^2 \right).
\end{equation}

We also introduce 
\begin{eqnarray*}
\Lambda_6(t)&=& j\left(\mathbf{Z}(t), \nabla^h y_h(\mathbf{Z}(t))\right)  
\frac{\mathbf{Z}^\prime(t)}{|\mathbf{Z}^\prime(t)|} \in \mathbb{R}^2\\
\Lambda_7(t)&=& \left( y_n(\mathbf{Z}(t)) \right)^2 
\frac{\mathbf{Z}^\prime(t)}{|\mathbf{Z}^\prime(t)|} \in \mathbb{R}^2 
\end{eqnarray*}
and the vectors:\\
$\widetilde{\Lambda}_6=\left(\widetilde{\Lambda}_6^1,\widetilde{\Lambda}_6^2\right)
\in \mathbb{R}^m \times \mathbb{R}^m$ with the components 
$\Lambda_6 (t_{k})-\Lambda_6 (t_{k+1})$, $1\leq k\leq m-1$ and the last component $\Lambda_6 (t_{m})$\\
$\widetilde{\Lambda}_7=\left(\widetilde{\Lambda}_7^1,\widetilde{\Lambda}_7^2\right)
\in \mathbb{R}^m \times \mathbb{R}^m$ with the components 
$\Lambda_7 (t_{k})-\Lambda_7 (t_{k+1})$, $1\leq k\leq m-1$ and the last component $\Lambda_7 (t_{m})$.
The last line of (\ref{3.14}) is approached by
\begin{equation}\label{4.16}
(\widetilde{\Lambda}_6^1)^T W^1 + (\widetilde{\Lambda}_6^2)^T W^2  
+\frac{1}{\epsilon}\left( (\widetilde{\Lambda}_7^1)^T W^1 + (\widetilde{\Lambda}_7^2)^T W^2 \right).
\end{equation}

\begin{proposition}\label{prop:4.1}
The discrete version of the relation (\ref{3.14}) is 
\begin{eqnarray}
\qquad dJ_{(G,U)}(R,V) & = &  
(\widetilde{\Lambda}_1^1)^T W^1 + (\widetilde{\Lambda}_1^2)^T W^2  
+(\widetilde{\Lambda}_2^1)^T W^1 + (\widetilde{\Lambda}_2^2)^T W^2
\nonumber\\
&+&(F^3)^T\ Q
+\frac{2}{\epsilon}\left( (\widetilde{\Lambda}_4^1)^T W^1 + (\widetilde{\Lambda}_4^2)^T W^2 \right)
+\frac{2}{\epsilon}Y^T N(Z)Q
\nonumber\\
&+& 
(\widetilde{\Lambda}_6^1)^T W^1 + (\widetilde{\Lambda}_6^2)^T W^2  
+\frac{1}{\epsilon}\left( (\widetilde{\Lambda}_7^1)^T W^1 + (\widetilde{\Lambda}_7^2)^T W^2 \right).
\label{4.17}
\end{eqnarray}
\end{proposition}

\noindent
\textbf{Proof.}
We get (\ref{4.17}) by adding (\ref{4.10}), (\ref{4.11}), (\ref{4.15}) and (\ref{4.16}).
\quad$\Box$

\bigskip

We point out that $Q$ depends on $V,R$ and $W$ depends on $R$, but $\widetilde{\Lambda}_i$,
$F^3$, $N(Z)$ as well as $Y$ are independent of $V,R$.

From (\ref{4.9}), we get
\begin{equation}\label{4.18}
Q=K^{-1}B^1(G)V + K^{-1}C^1(G,U)R
\end{equation}
and the discrete version of the operator $A $ in the Corollary \ref{cor:3.6} is
$$
(R,V)\in \mathbb{R}^n \times \mathbb{R}^n \rightarrow
A^1(R,V)=K^{-1}B^1(G)V + K^{-1}C^1(G,U)R.
$$

Next, we present how $W$ depends on $R$.
Let us introduce the square matrices of order 2
$$
A_2(k)=\left(
\begin{array}{cc}
-\Delta t\, \partial_1^h\partial_2^h g_h(Z_{k+1}) & -\Delta t\, \partial_2^h\partial_2^h g_h(Z_{k+1})\\
\Delta t\, \partial_1^h\partial_1^h g_h(Z_{k+1}) & \Delta t\, \partial_2^h\partial_1^h g_h(Z_{k+1})
\end{array}
\right),
$$
$$
I_2=\left(
\begin{array}{rr}
  1 & 0\\
  0 & 1
\end{array}
\right),\qquad
M_2(k)=\left( I_2 -\Delta tA_2(k)\right)^{-1}
$$
and the $2\times n$ matrice
$$
N_2(k) = \left(
\begin{array}{r}
-\Delta t\, \Phi^T(Z_{k+1}) \Pi_h^2\\
\Delta t\, \Phi^T(Z_{k+1}) \Pi_h^1
\end{array}
\right)
$$
where $\Phi(Z_k)=(\phi_i(Z_k))_{i\in I}^T \in \mathbb{R}^n$.
We can rewrite the system (\ref{4.12})-(\ref{4.13}) as 
$$
\left(
\begin{array}{c}
W_{k+1}^1 \\
W_{k+1}^2
\end{array}
\right)
= M_2(k)
\left(
\begin{array}{c}
W_{k}^1 \\
W_{k}^2
\end{array}
\right)
+M_2(k) N_2(k)R.
$$
We have the following equality
\begin{eqnarray}
\left(\begin{array}{c}
\begin{array}{c}
W_{1}^1\\
W_{1}^2
\end{array}
\\
\vdots
\\
\begin{array}{c}
W_{m}^1 \\
W_{m}^2
\end{array}
\end{array}
\right)
& = & M_{2m}
\times
\left(
\begin{array}{c}
N_2(0) \\
N_2(1) \\
\vdots  \\
N_2(m-1)
\end{array}
\right)
R
\label{4.19}
\end{eqnarray}
the right-hand side, $M_{2m}$ is a square matrix of order $2m$ given by
$$
\left(
\begin{array}{ccccc}
M_2(0)       & 0      & \cdots & 0 & 0\\
M_2(1)M_2(0) & M_2(1) & \cdots & 0 & 0\\
\vdots & \vdots & & \vdots & \vdots \\
\Pi_{k=0}^{m-1} M_2(k), & \Pi_{k=1}^{m-1} M_2(k) , & \cdots & \Pi_{k=m-2}^{m-1} M_2(k), & M_2(m-1)
\end{array}
\right)
$$
and the second matrix, which contains $N_2$, is of size $2m\times n$.
Now, $W$ depends on $R$ by (\ref{4.19}), we define
the linear operator approximation of $B$ from the Corollary \ref{cor:3.6}
\begin{equation}\label{4.20}
R \in \mathbb{R}^n \rightarrow W=\left( W^1, W^2 \right) = \left( B^2(G,Z) R,
B^3(G,Z) R \right) \in \mathbb{R}^m \times \mathbb{R}^m.
\end{equation}

We can rewrite (\ref{4.17}) as
\begin{eqnarray}
\qquad dJ_{(G,U)}(R,V) & = &  
(\widetilde{\Lambda}_1^1 + \widetilde{\Lambda}_2^1 + \frac{2}{\epsilon}\widetilde{\Lambda}_4^1)^T 
B^2(G,Z)R
\nonumber\\
&+&
(\widetilde{\Lambda}_1^2 + \widetilde{\Lambda}_2^2 + \frac{2}{\epsilon}\widetilde{\Lambda}_4^2)^T
B^3(G,Z)R
\nonumber\\
&+&(F^3+\frac{2}{\epsilon}Y^T N(Z))^T K^{-1}B^1(G)V
\nonumber\\
&+&(F^3+\frac{2}{\epsilon}Y^T N(Z))^T K^{-1}C^1(G,U)R
\nonumber\\
&+& 
(\widetilde{\Lambda}_6^1 +\frac{1}{\epsilon} \widetilde{\Lambda}_7^1)^T B^2(G,Z)R
\nonumber\\
&+&(\widetilde{\Lambda}_6^2 +\frac{1}{\epsilon} \widetilde{\Lambda}_7^2)^T B^3(G,Z)R
\label{4.21}
\end{eqnarray}

The first four lines of (\ref{4.21}) represent an approximation of the
first four lines of (\ref{3.17}).


\textbf{Descent direction method}

The descent direction method needs at each step a descent direction, i.e.
$(R^k,V^k)$ such that $dJ_{(G^k,U^k)}(R^k,V^k) <0$ and the next step is defined by
$$
(G^{k+1},U^{k+1})=(G^k,U^k)+\lambda_k (R^k,V^k),
$$
where $\lambda_k > 0$ is computed by some line search
$$
\lambda_k \in \arg\min_{\lambda >0} J \left( (G^k, U^k) + \lambda (R^k, V^k) \right) .
$$
The algorithm stops if $dJ_{(G^k,U^k)}(R^k,V^k)=0$ or
$| J(G^{k+1},U^{k+1}) - J(G^k,U^k)| < tol$ for some prescribed 
tolerance parameter $tol$.

\begin{proposition}\label{prop:4.2}
A descent direction for $J$ at $(G,U)$ is 
$(R^{*},V^*)\in \mathbb{R}^n \times \mathbb{R}^n$ given by
\begin{eqnarray*}
(V^*)^T&=&
-(F^3+\frac{2}{\epsilon}Y^T N(Z))^T K^{-1}B^1(G)
\nonumber\\
(R^{*})^T&=&
-(\widetilde{\Lambda}_1^1 + \widetilde{\Lambda}_2^1 + \frac{2}{\epsilon}\widetilde{\Lambda}_4^1)^T 
B^2(G,Z)
\nonumber\\
&-&
(\widetilde{\Lambda}_1^2 + \widetilde{\Lambda}_2^2 + \frac{2}{\epsilon}\widetilde{\Lambda}_4^2)^T
B^3(G,Z)
\nonumber\\
&-&(F^3+\frac{2}{\epsilon}Y^T N(Z))^T K^{-1}C^1(G,U)
\nonumber\\
&-& 
(\widetilde{\Lambda}_6^1 +\frac{1}{\epsilon} \widetilde{\Lambda}_7^1)^T B^2(G,Z)
\nonumber\\
&-&(\widetilde{\Lambda}_6^2 +\frac{1}{\epsilon} \widetilde{\Lambda}_7^2)^T B^3(G,Z) . 
\end{eqnarray*}
\end{proposition}  

\noindent
\textbf{Proof.}
We can rewrite (\ref{4.21}) as
$dJ_{(G,U)}(R,V)=-(V^*)^T V - (R^{*})^T R$, then\\
$dJ_{(G,U)}(R^{*},V^*)=-\| V^*\|_{\mathbb{R}^n}^2 -\| R^{*}\|_{\mathbb{R}^n}^2 \leq 0$.
If the gradient  $dJ_{(G,U)} = (R^{*},V^*)$ is non null (non stationary points),
the inequality is strict.
\quad$\Box$  

\bigskip

Let us introduce a simplified adjoint system: find $p_h$ in $\mathbb{V}_h$ such that
\begin{eqnarray}
\int_D \nabla \varphi_h \cdot \nabla p_h d\mathbf{x}
&= &\int_0^{T_g} \nabla_2 j\left(Z(t), y_h(Z(t))\right)\cdot
\nabla^h \varphi_h(Z(t))|Z^\prime(t)| dt
\nonumber\\
&+&\frac{2}{\epsilon}\int_0^{T_g} y_h(Z(t))\varphi_h(Z(t))|Z^\prime(t)| dt
\label{4.22}
\end{eqnarray}
$\forall \varphi_h \in \mathbb{V}_h$ and $Z(t)$ satisfying (\ref{4.3})-(\ref{4.5}).
We have $p_h=\sum_{i\in I_0} P_i \phi_i$ and
$P=(P_i)_{i\in I_0}^T\in \mathbb{R}^{n_0}$.

\begin{proposition}\label{prop:4.3}
Given $g_h,u_h\in\mathbb{W}_h$, let $y_h\in\mathbb{V}_h$ be the solution of
(\ref{4.1}). For $r_h=-p_hu_h$, $v_h=-p_h$, with $p_h\in\mathbb{V}_h$
the solution of
(\ref{4.22}), then
\begin{eqnarray}
&&\int_0^{T_g} \nabla_2 j\left(Z(t), y_h(Z(t))\right)\cdot \nabla^h q_h(Z(t))|Z^\prime(t)| dt
\nonumber\\ 
&+&\frac{2}{\epsilon}\int_0^{T_g} y_h(Z(t))q_h(Z(t))|Z^\prime(t)| dt
\leq 0,
\label{4.23}
\end{eqnarray}
where $q_h \in \mathbb{V}_h$ is the solution of (\ref{4.8}) depending on $r_h $
and $v_h $.
\end{proposition}

\noindent\textbf{Proof.}
Putting $\varphi_h=p_h$ in (\ref{4.8}) and $\varphi_h=q_h$ in (\ref{4.22}), we get
\begin{eqnarray*}
&&\int_D \left( (g_h)_+^2 v_h + 2(g_h)_+ u_h r_h\right) p_h d\mathbf{x}  
 = \int_D \nabla q_h \cdot \nabla p_h d\mathbf{x}\\
&=&  \int_0^{T_g} \nabla_2 j\left(Z(t), y_h(Z(t))\right)\cdot
\nabla^h q_h(Z(t))|Z^\prime(t)| dt\\
&+&\frac{2}{\epsilon}\int_0^{T_g} y_h(Z(t))q_h(Z(t))|Z^\prime(t)| dt .
\end{eqnarray*}

For $v_h=-p_h$, we have
$$
\int_D (g_h)_+^2 v_h p_h d\mathbf{x}
=-\int_D (g_h)_+^2 p_h^2d\mathbf{x}
\leq 0
$$
and for $r_h=-p_h\,u_h$, we have
$$
\int_D  2(g_h)_+ u_h r_h p_h d\mathbf{x}
=-\int_D  2(g_h)_+ (u_h p_h)^2 d\mathbf{x}
\leq 0
$$
since $(g_h)_+ \geq 0$ in $D$.
\quad$\Box$

\medskip

\begin{remark}\label{rem:4.1}

The terms from (\ref{3.14}), (\ref{4.22}), (\ref{4.23}), containing $q$, can be rewritten as integrals over $\partial \Omega_g$. For instance, for (\ref{4.22}) we obtain:

\begin{eqnarray}
\int_D \nabla \varphi_h \cdot \nabla p_h d\mathbf{x}
&= &\int_{\partial \Omega_{g_h}} \nabla_2 j\left(s, y_h(s)\right)\cdot
\nabla^h \varphi_h(s) ds
\nonumber\\
&+&\frac{2}{\epsilon}\int_{\partial \Omega_{g_h}} y_h(s)\varphi_h(s) ds.
\nonumber
\end{eqnarray}

However, the way they are expressed in (\ref{3.14}), (\ref{4.22}), (\ref{4.23}) avoids the use of the unknown geometry and all the elements are easily computable. For instance, $T_g$ is obtained automatically when solving the Hamiltonian system (\ref{2.10})-(\ref{2.12}), while the initial conditions (on each component of $G$) are simply obtained via the equation $g = 0$ and standard routines, together with a simple iterative procedure to generate all of them. See as well \cite{MT2019}.

\end{remark}

\section{Numerical tests}
\setcounter{equation}{0}

In the numerical examples, we have employed the software FreeFem++, \cite{freefem++}.

The functional appearing in the objective function is
$$
j\left(\mathbf{x}, \nabla y(\mathbf{x})\right)
=\frac{1}{2}
\left(
\frac{\partial y}{\partial \mathbf{n}}(\mathbf{x}) - \delta(\mathbf{x})
\right)^2
$$
where $\delta\in H^1(D)$ is a given function.
It follows that
\begin{eqnarray*}
  \nabla_1 j\left(\mathbf{x}, \nabla y(\mathbf{x})\right)
  &=&-\left(
\frac{\partial y}{\partial \mathbf{n}}(\mathbf{x}) - \delta(\mathbf{x})
\right)\nabla \delta(\mathbf{x})\\
\nabla_2 j\left(\mathbf{x}, \nabla y(\mathbf{x})\right)
&=&\left(
\frac{\partial y}{\partial \mathbf{n}}(\mathbf{x}) - \delta(\mathbf{x})
\right)\mathbf{n}.
\end{eqnarray*}

\newpage

\textbf{Example 1.}

a) The computational domain is $D=]-1,1[\times ]-1,1[$,
the load is $f=-4$ and $\delta=1$.
This problem has the solution $y_e(x_1,x_2)=x_1^2 + x_2^2 -0.5^2$
defined on the disk of center $(0,0)$ and radius $0.5$.
The mesh of $D$ has 53290 triangles and 26946 vertices.
The penalization parameter is $\epsilon=10^{-4}$ and the tolerance
parameter for the stopping test is $ tol = 10^{-6} $.
The initial domain is the disk of center $(0.2,0.2)$ and radius $0.5$,
given by
\begin{equation}\label{5.1}
g_0(x_1,x_2)=(x_1-0.2)^2 + (x_2-0.2)^2 -0.5^2.
\end{equation}

As descent direction, we use $(R^k,V^k)$ given by Proposition \ref{prop:4.3}.
For $r_h$, $v_h$ given by Proposition \ref{prop:4.3}
and a scaling parameter $\gamma >0$,
then $\gamma r_h$ and $v_h$ also give a descent direction.
We take here $\gamma=\frac{1}{\|r_h\|_\infty}$. 

\begin{figure}[ht]
\begin{center}
\includegraphics[width=4.5cm]{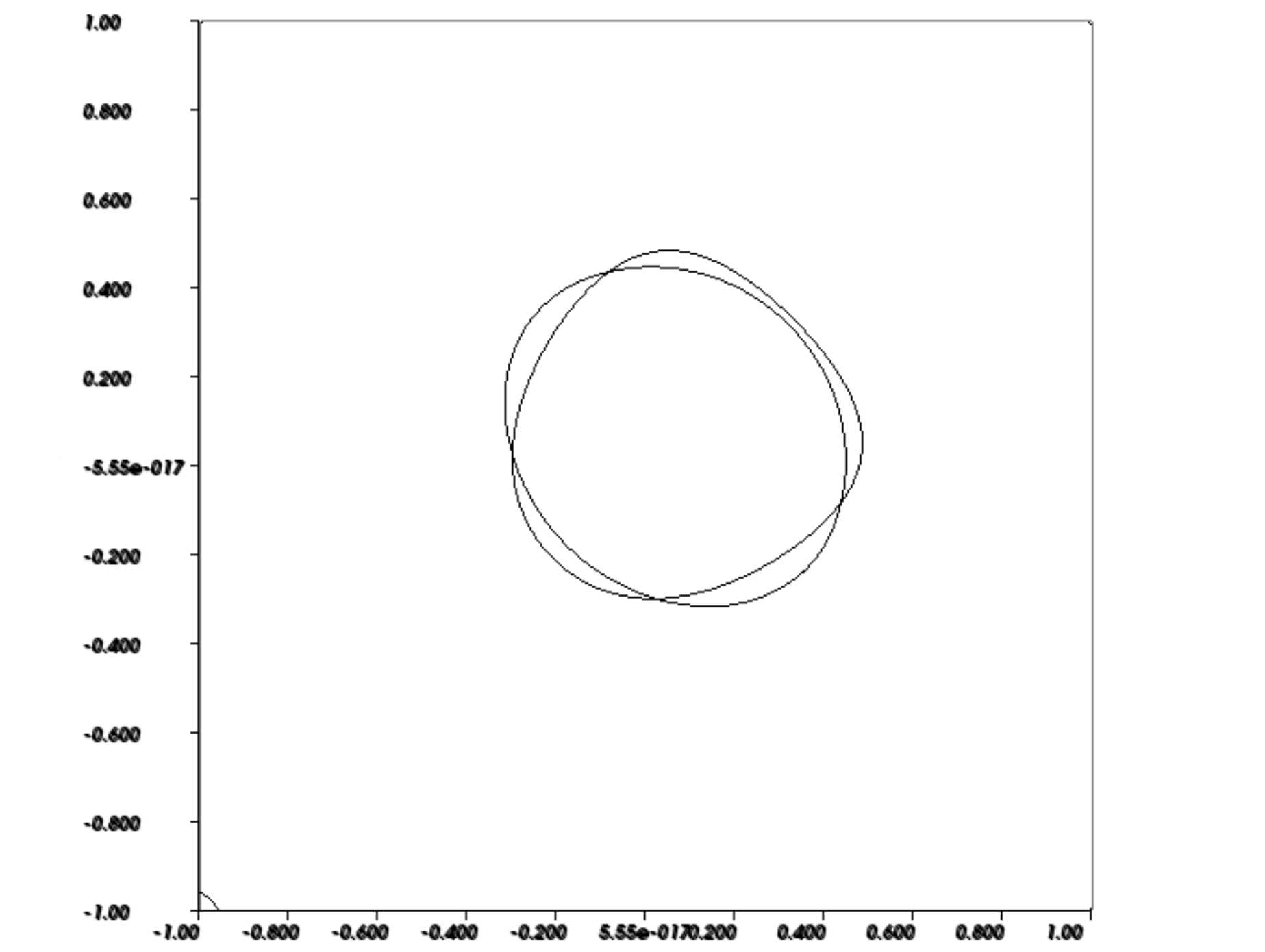}
\   
\includegraphics[width=6.5cm]{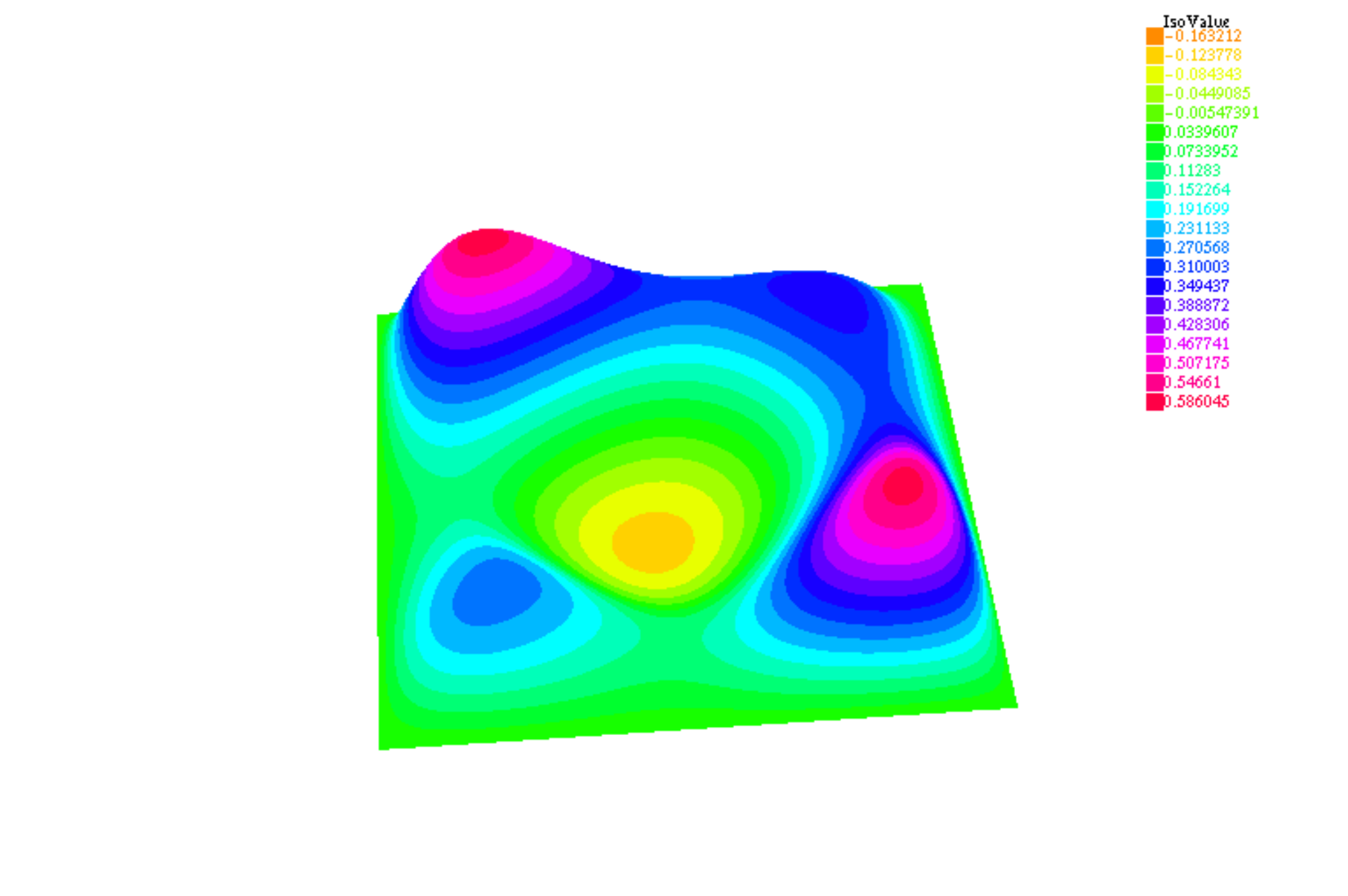}
\
\includegraphics[width=5.5cm]{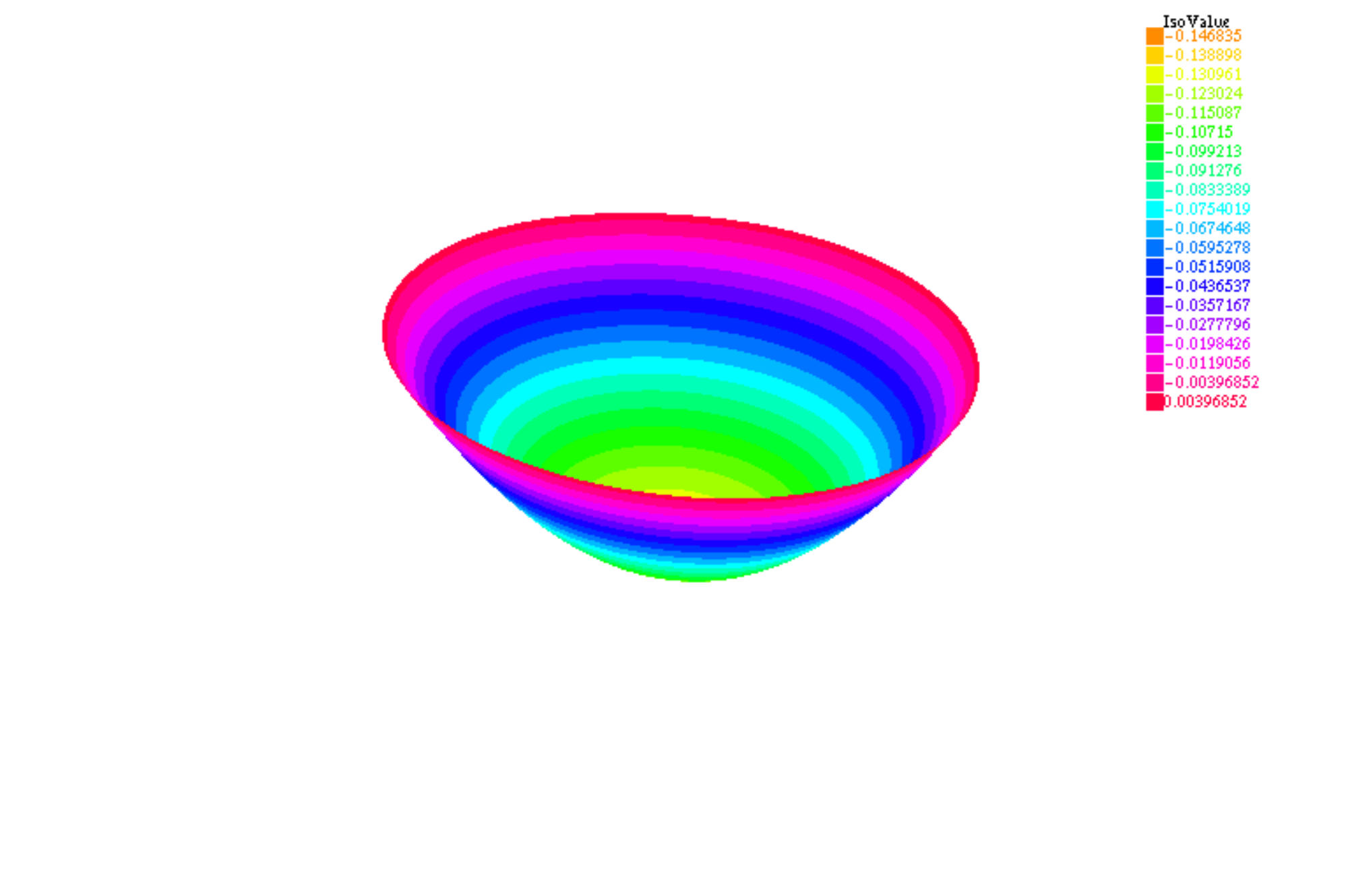}
\
\includegraphics[width=5.5cm]{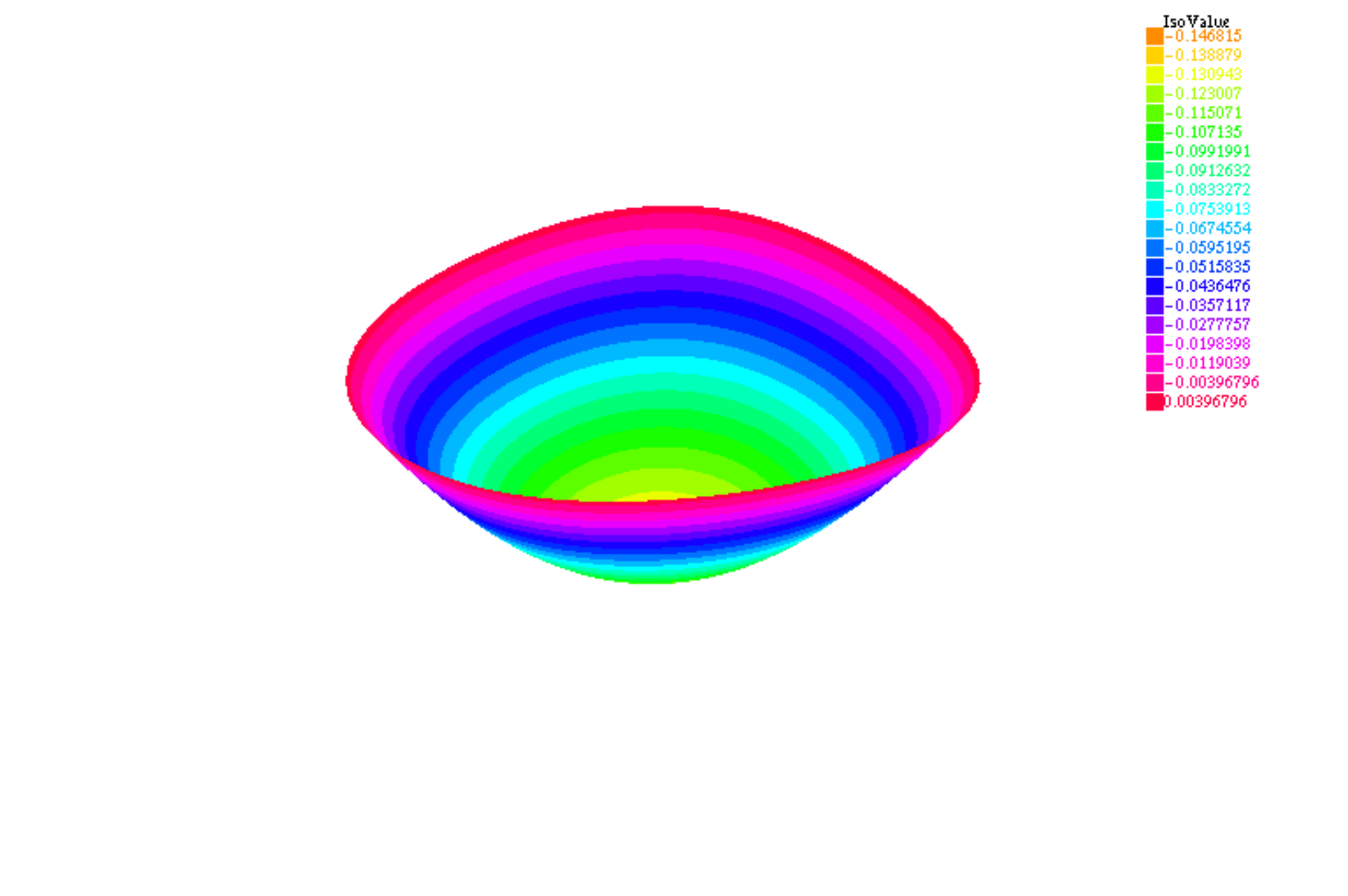}
\end{center}
\caption{Example 1a. 
The zero level sets of the computed optimal $g$, $y$ (top, left),
the final state $y$ (top, right),
the solution of the system (\ref{1.2})-(\ref{1.3})
in $\Omega_g$ (bottom, left) and
in the domain of boundary the zero level sets of $y$ (bottom, right).
\label{fig:ex1_k_15}}
\end{figure}

\noindent
Notice that the difference between the two curves (Figure \ref{fig:ex1_k_15}, top, left)
is due to the fact that the penalization integral is not null at the final step.

\begin{figure}[ht]
\begin{center}
\includegraphics[width=4.5cm]{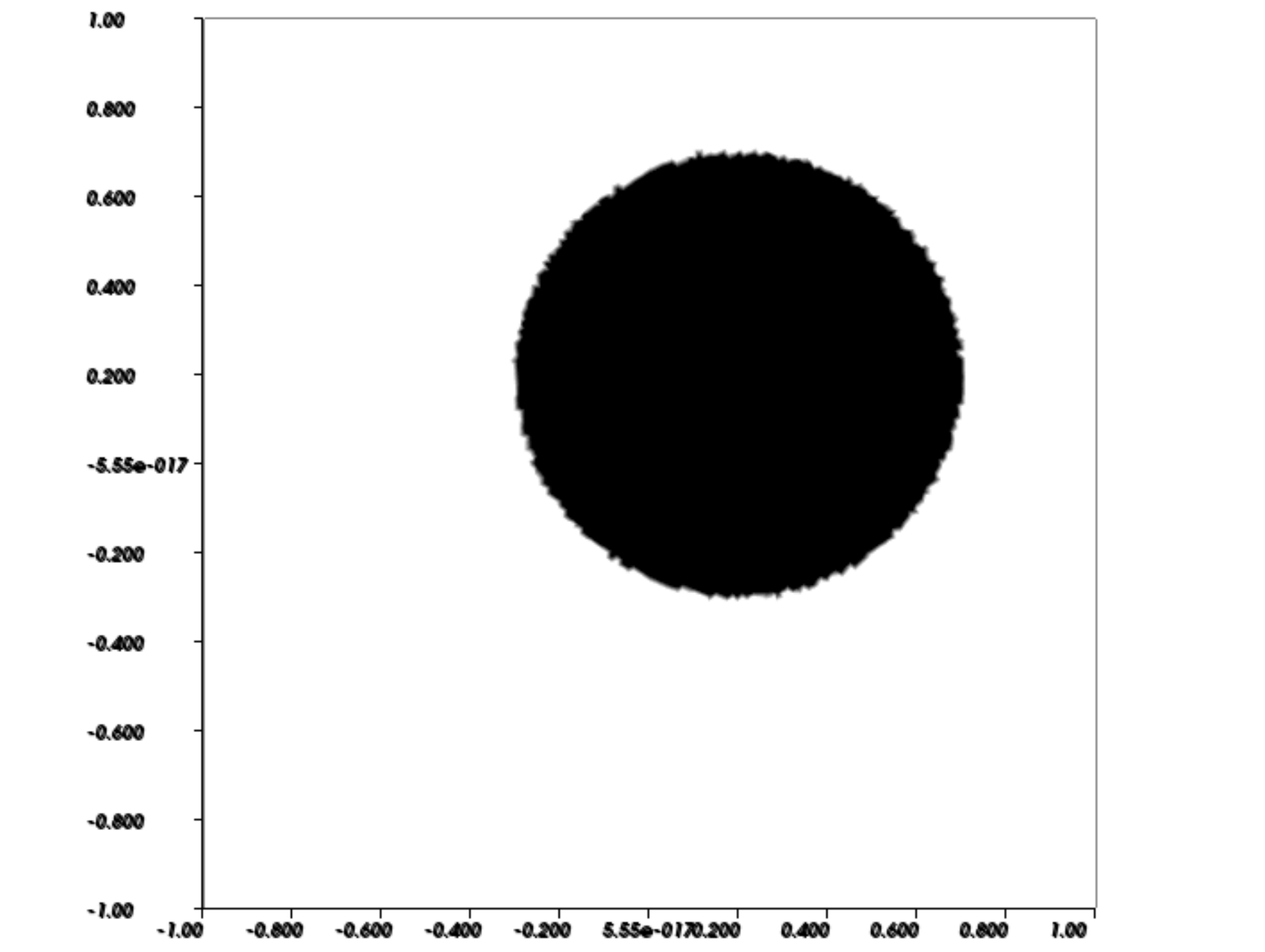}
\   
\includegraphics[width=4.5cm]{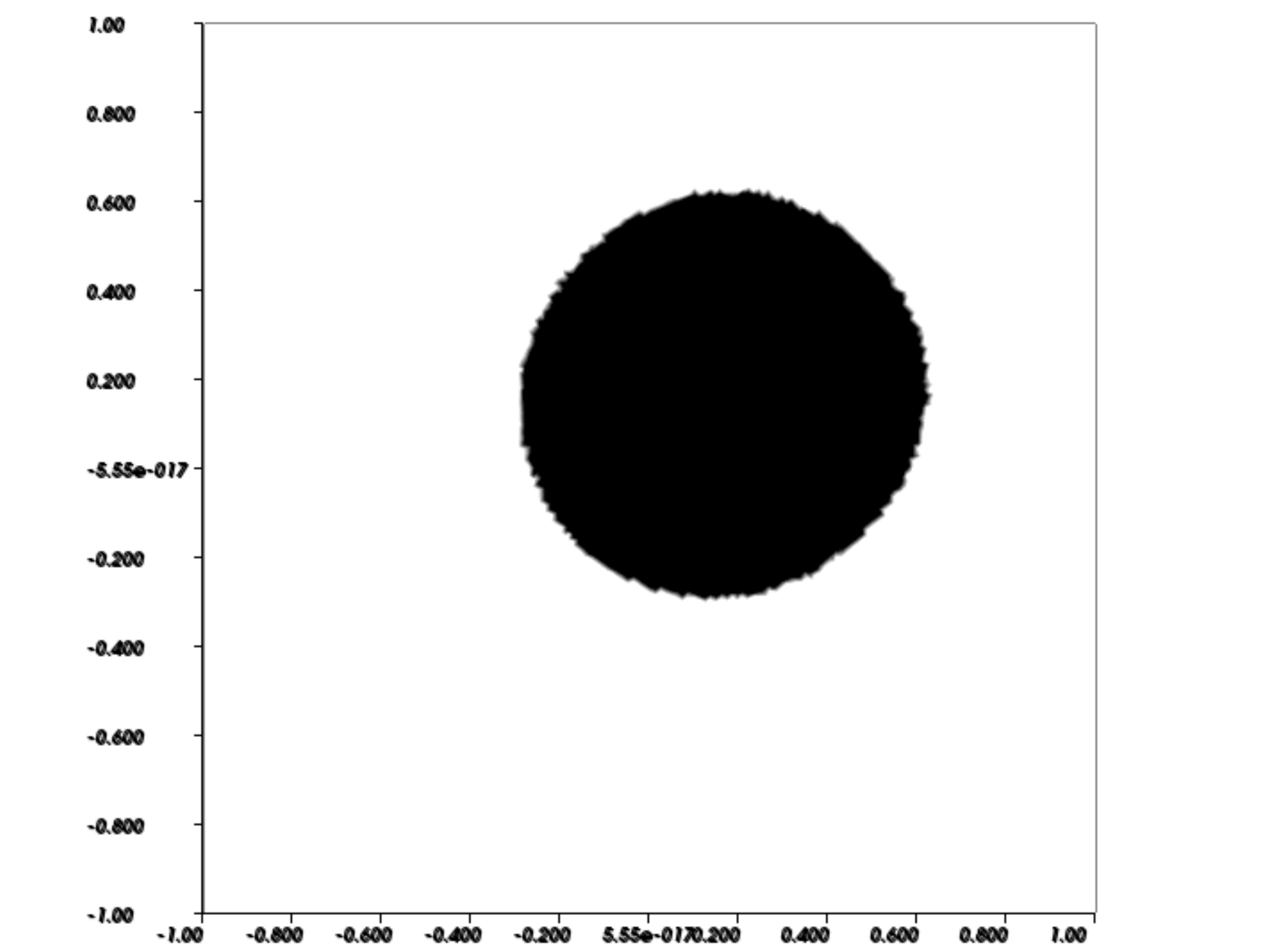}\\
\includegraphics[width=4.5cm]{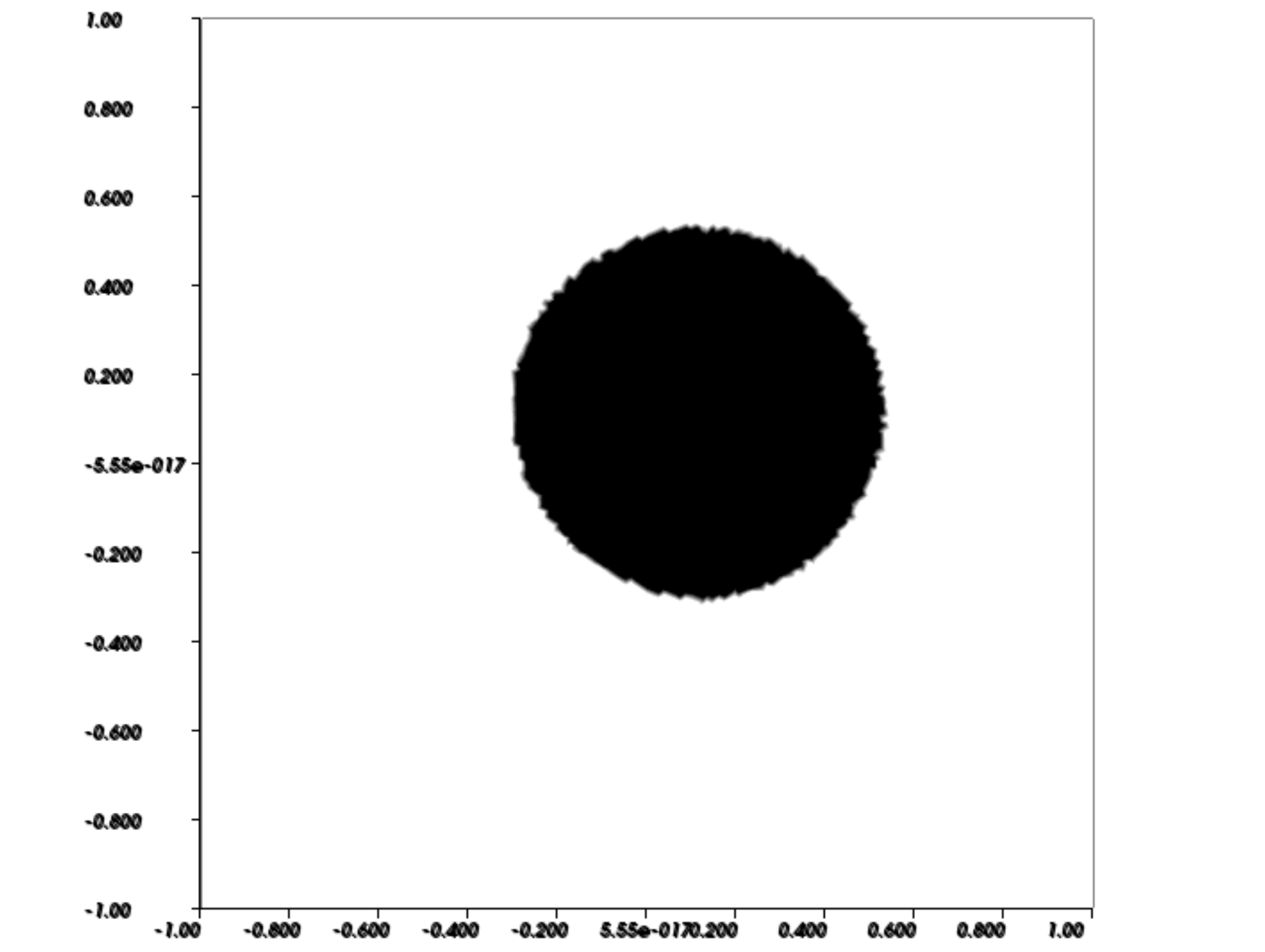}
\
\includegraphics[width=4.5cm]{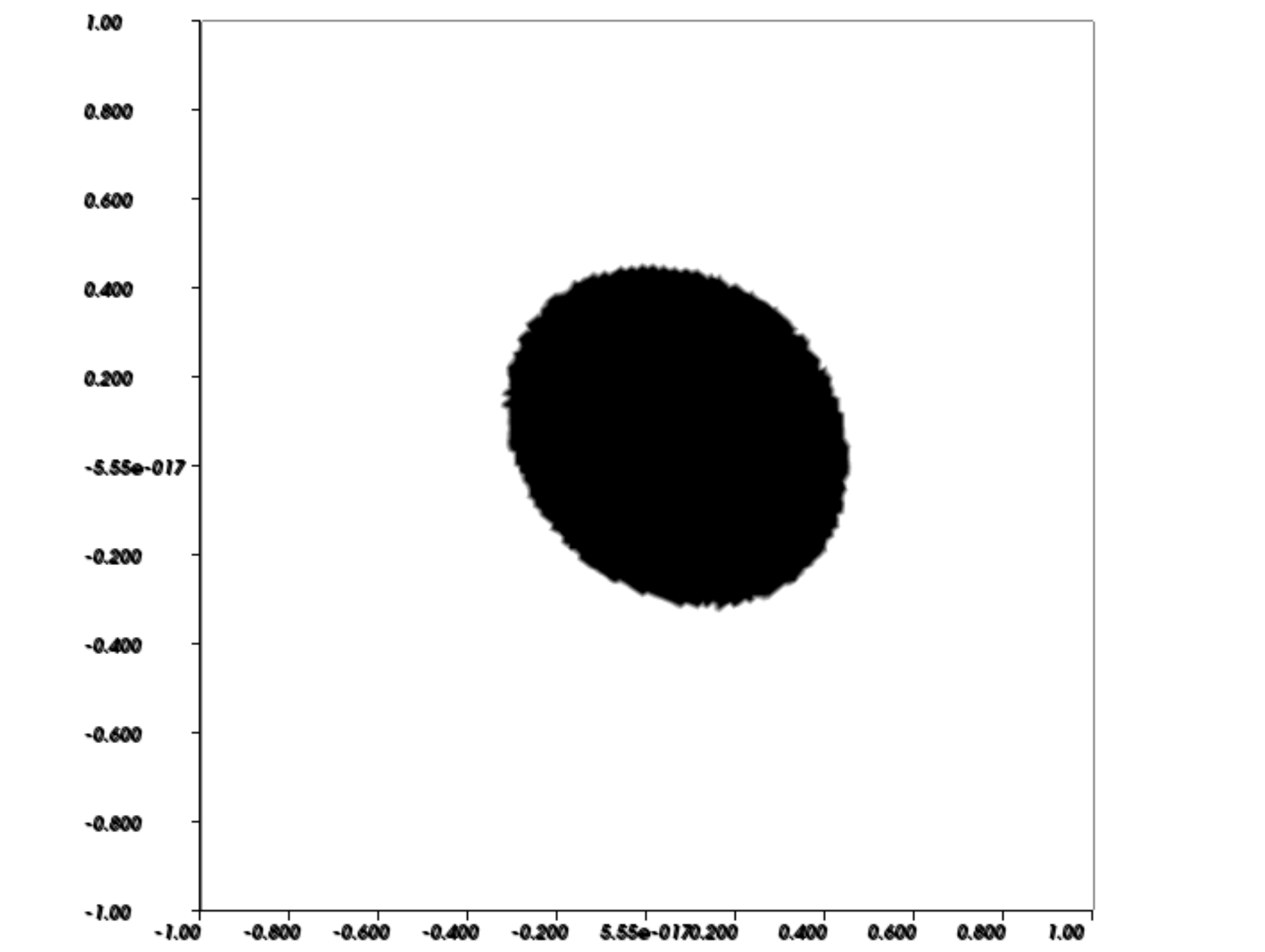}
\end{center}
\caption{Example 1a. Initial domain (top, left), for $k=4$ (top, right),
for $k=8$ (bottom, left) and the final domain (bottom, right).
\label{fig:ex1a_xi}}
\end{figure}

The stopping test is obtained for $k=13$.
The objective function (\ref{1.3}) is $0.072180$
for the solution of the elliptic system (\ref{1.2})-(\ref{1.3})
in the domain $\Omega_g$ and $0.077413$ for the solution in the domain of boundary
the zero level sets of $y$.
The initial, intermediate and the final domains are presented in Figure \ref{fig:ex1a_xi}
and the corresponding values of the objective function (\ref{3.1})
are detailed in Table \ref{tab:ex1a_J}.

\begin{table}[ht]
\begin{center}
\begin{tabular}{|c|c|c|c|c|}
\hline
iteration & k=0 & k=4 & k=8 & final \\  \hline
$t_1$   & 0.268404 & 0.259176 & 0.151267 & 0.093539 \\  \hline
$t_2$   & 2.38981 & 0.110803 & 0.036932 & 0.002103 \\  \hline
$J$   & 23898.4 & 1108.29 & 369.471 & 21.133 \\  \hline
\end{tabular}
\end{center}
\caption{Example 1a. The computed objective function (\ref{3.1}), i.e.
$J=t_1+\frac{1}{\epsilon}t_2$, where
$t_1=\int_{\partial\Omega_g}j\left(s, \nabla y(s)\right)ds$ and
$t_2=\int_{\partial\Omega_g}\left(y(s)\right)^2 ds$.}
\label{tab:ex1a_J}
\end{table}

\clearpage

\bigskip

b) We have the same parameters as before, just the 
initial domain is the disk of center $(0.2,0.2)$ and radius $0.4$
with a circular hole of center $(0.2,0.2)$ and radius $0.2$
with $g_0(x_1,x_2)$ given by
\begin{equation}\label{5.2}
\max\left(
 (x_1-0.2)^2 + (x_2-0.2)^2 -0.4^2,
-(x_1-0.2)^2 - (x_2-0.2)^2 +0.2^2 
\right) .
\end{equation}

\begin{figure}[ht]
\begin{center}
\includegraphics[width=4.5cm]{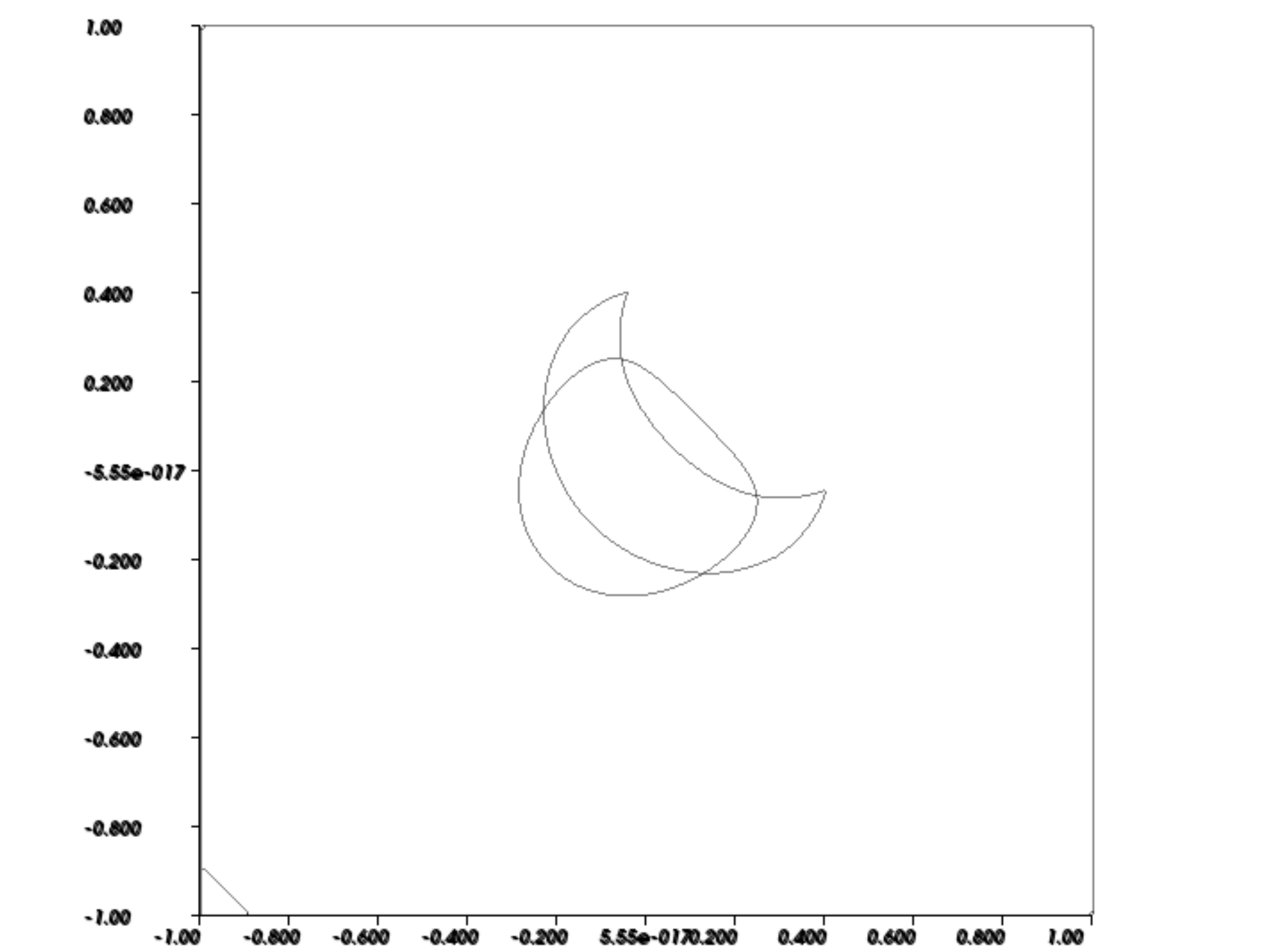}
\   
\includegraphics[width=6.5cm]{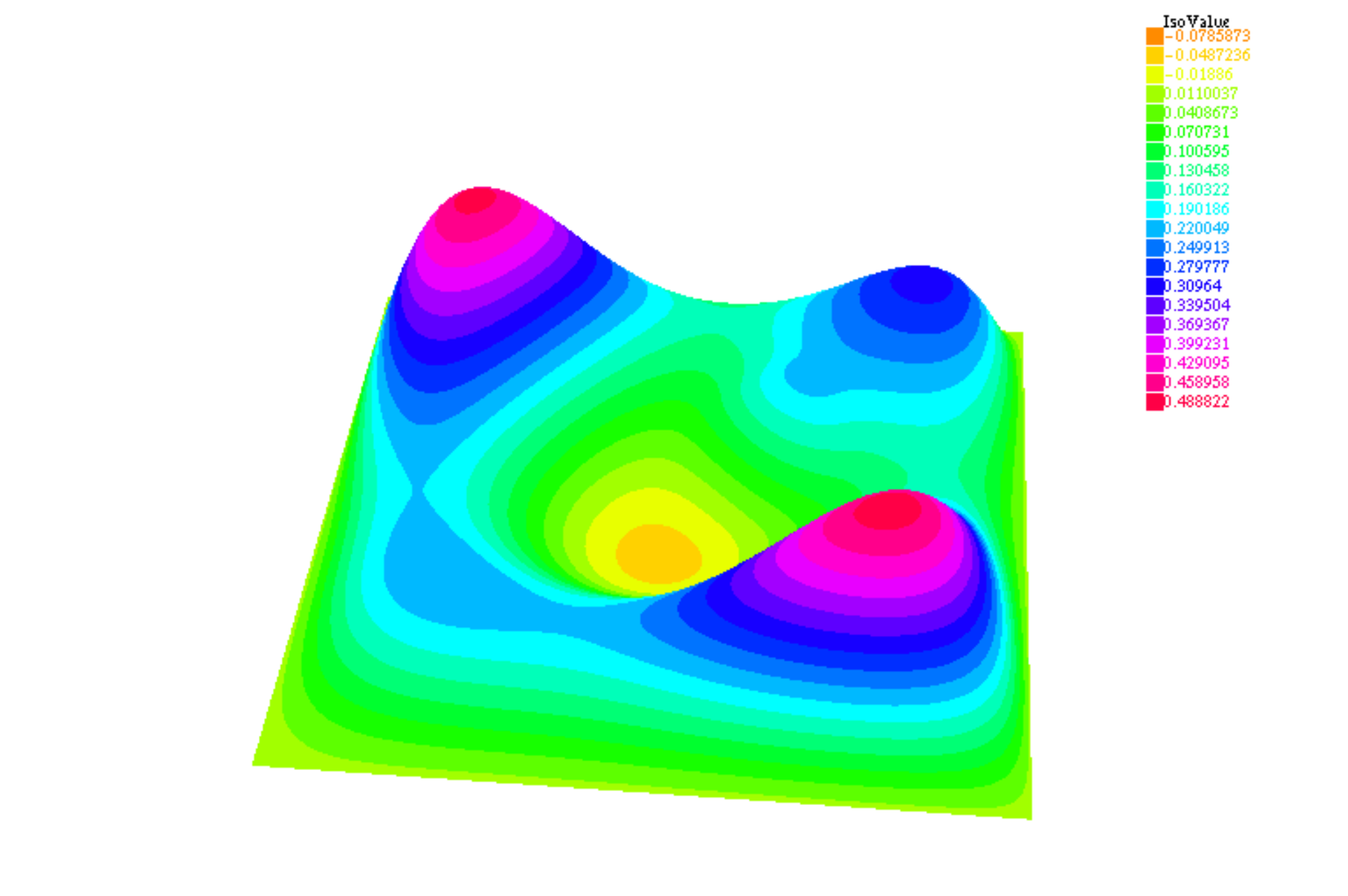}
\
\includegraphics[width=5.5cm]{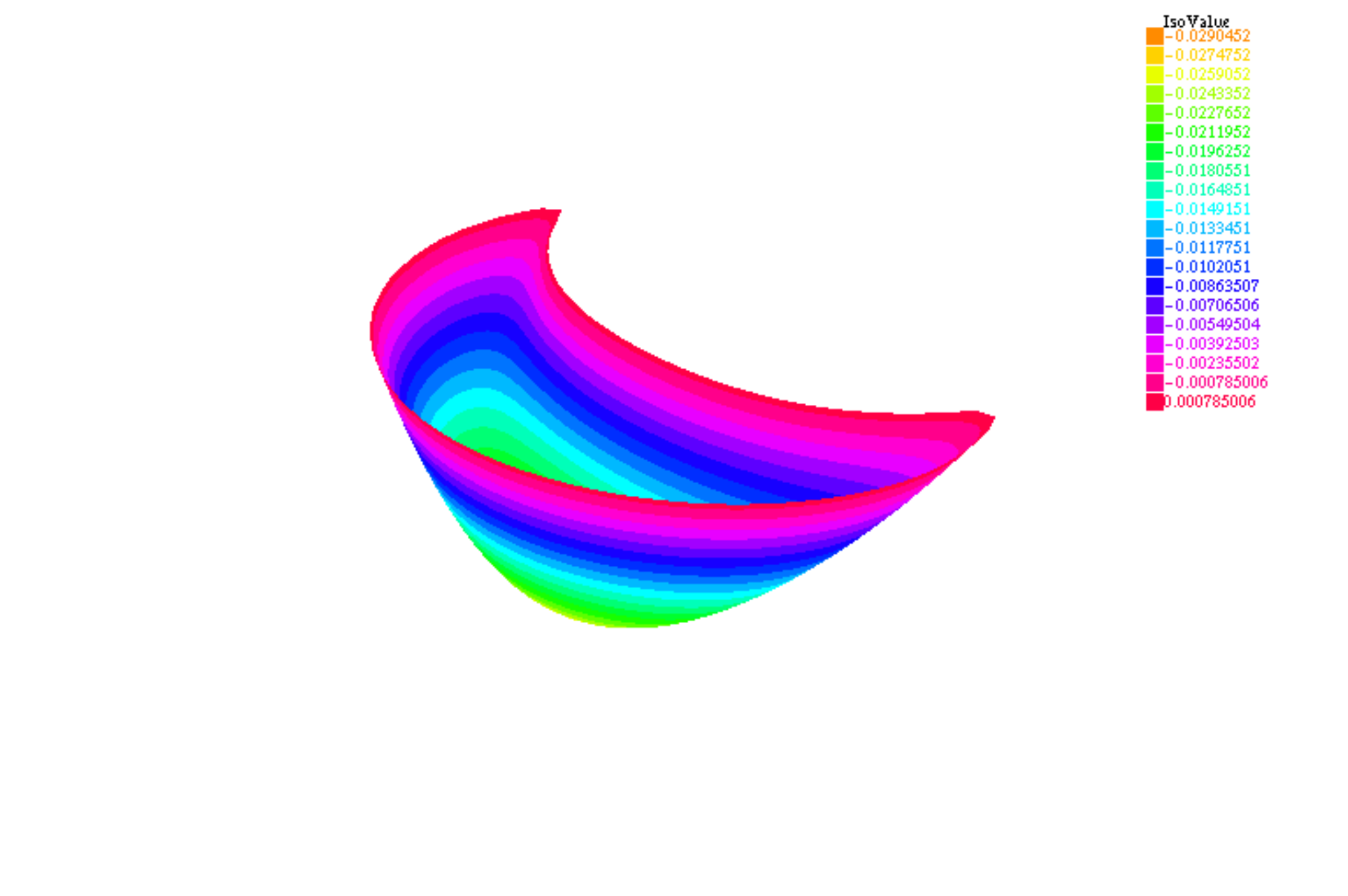}
\
\includegraphics[width=5.5cm]{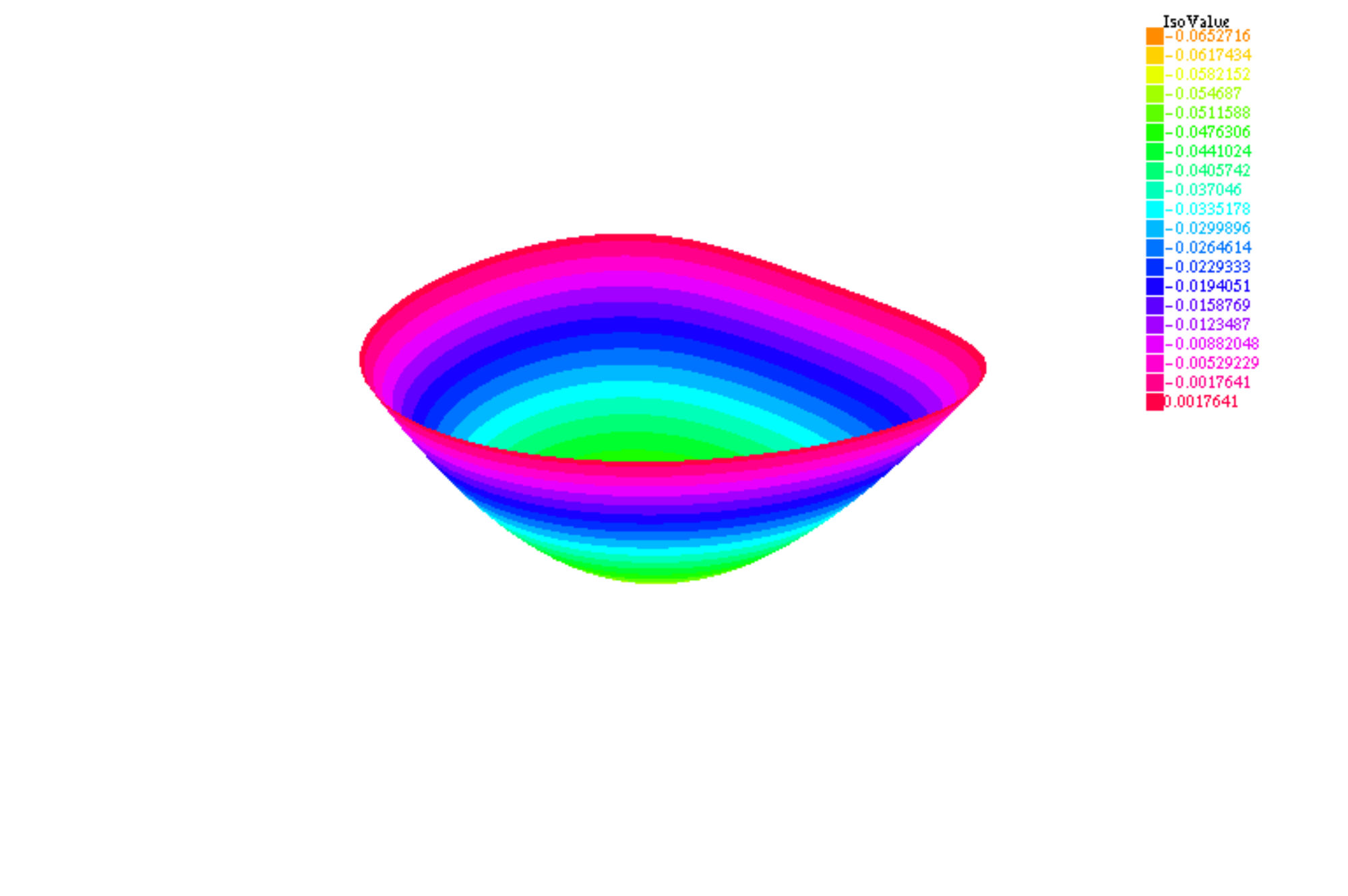}
\end{center}
\caption{Example 1b. 
The zero level sets of the computed optimal $g$, $y$ (top, left),
the final state $y$ (top, right),
the solution of the problem (\ref{1.2})-(\ref{1.3}) written
in $\Omega_g$ (bottom, left) and
in the domain of boundary the zero level sets of $y$ (bottom, right).
\label{fig:ex2_k_15}}
\end{figure}

\begin{figure}[ht]
\begin{center}
\includegraphics[width=4.5cm]{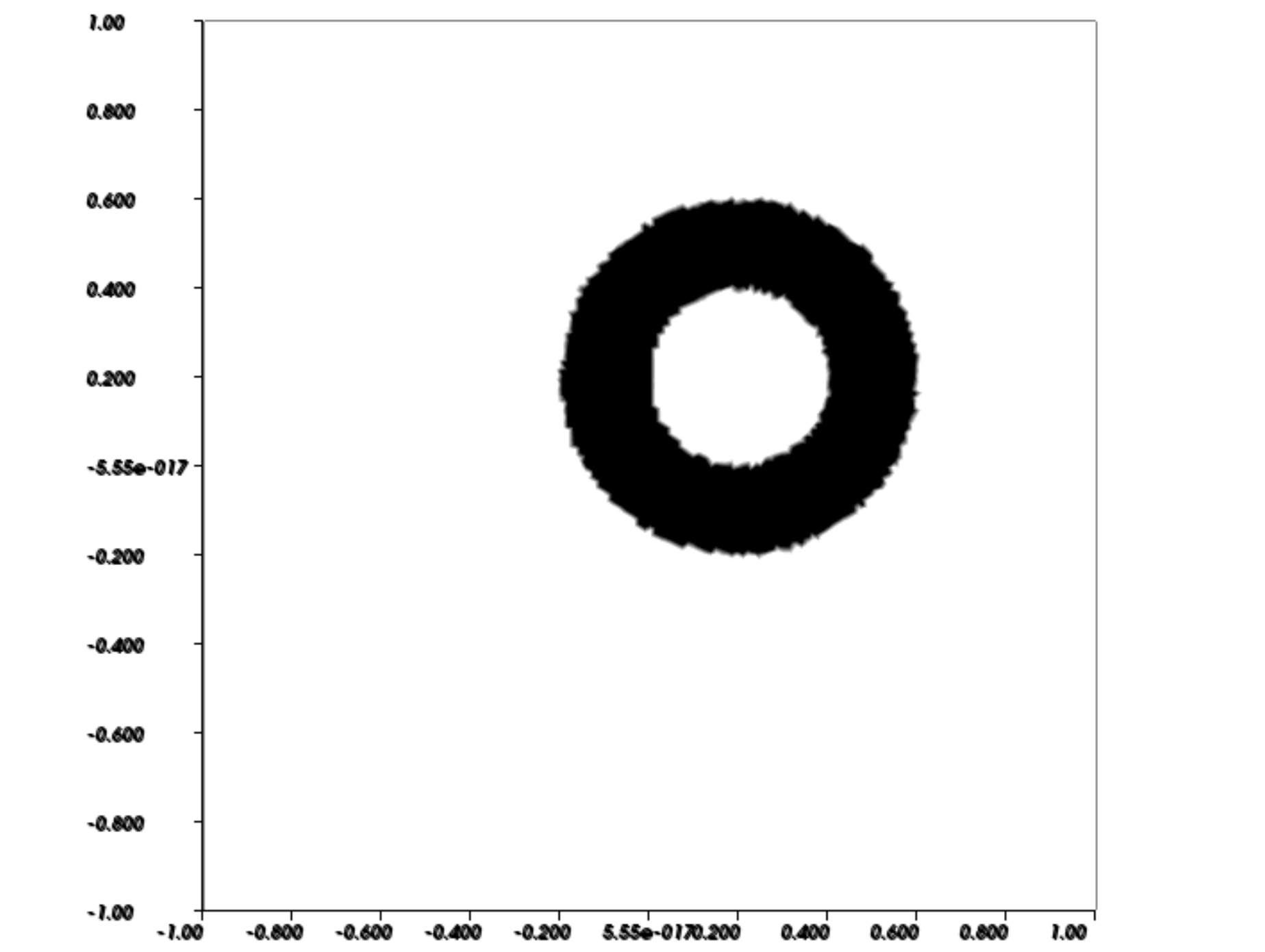}
\   
\includegraphics[width=4.5cm]{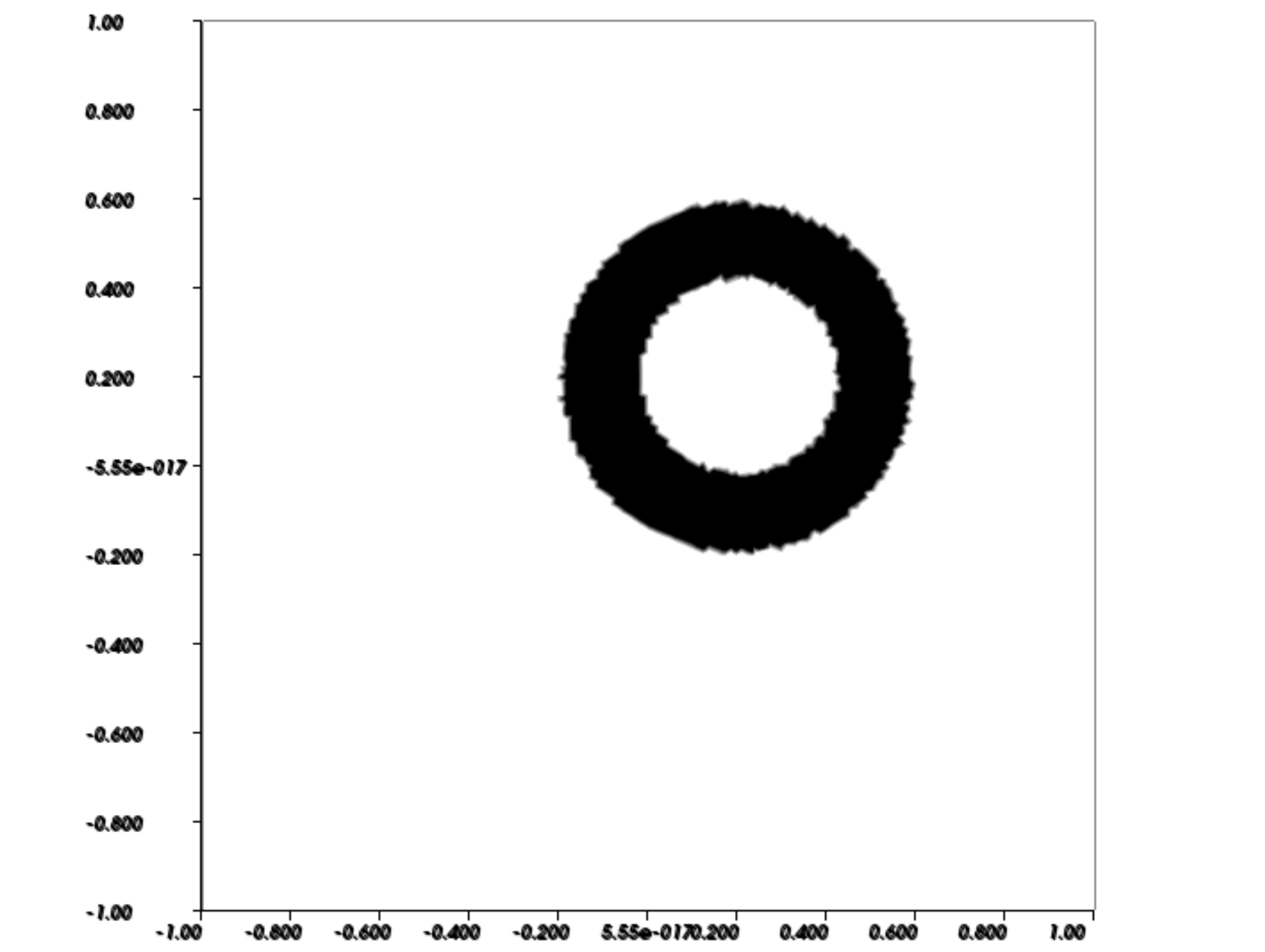}
\
\includegraphics[width=4.5cm]{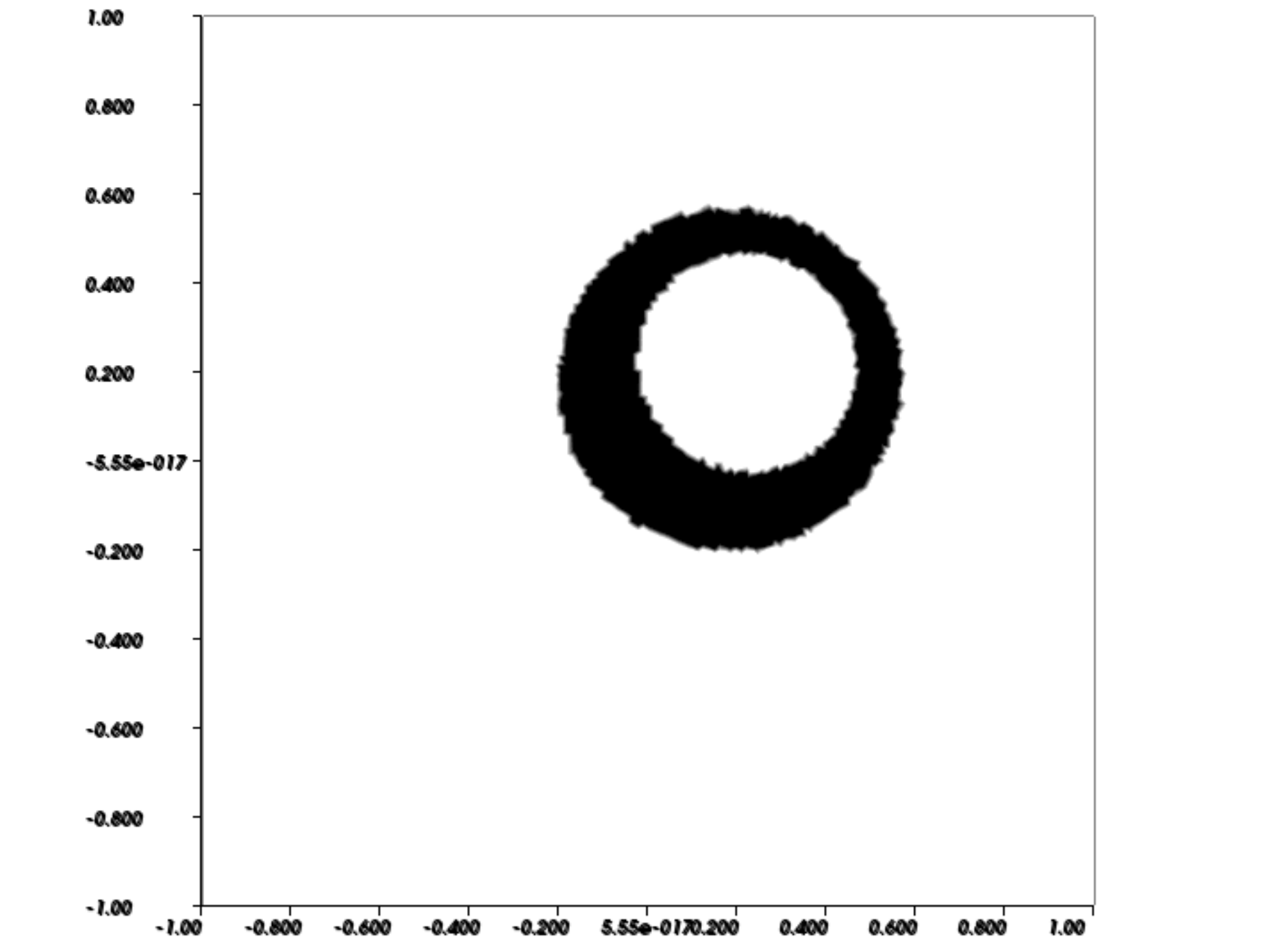}\\
\includegraphics[width=4.5cm]{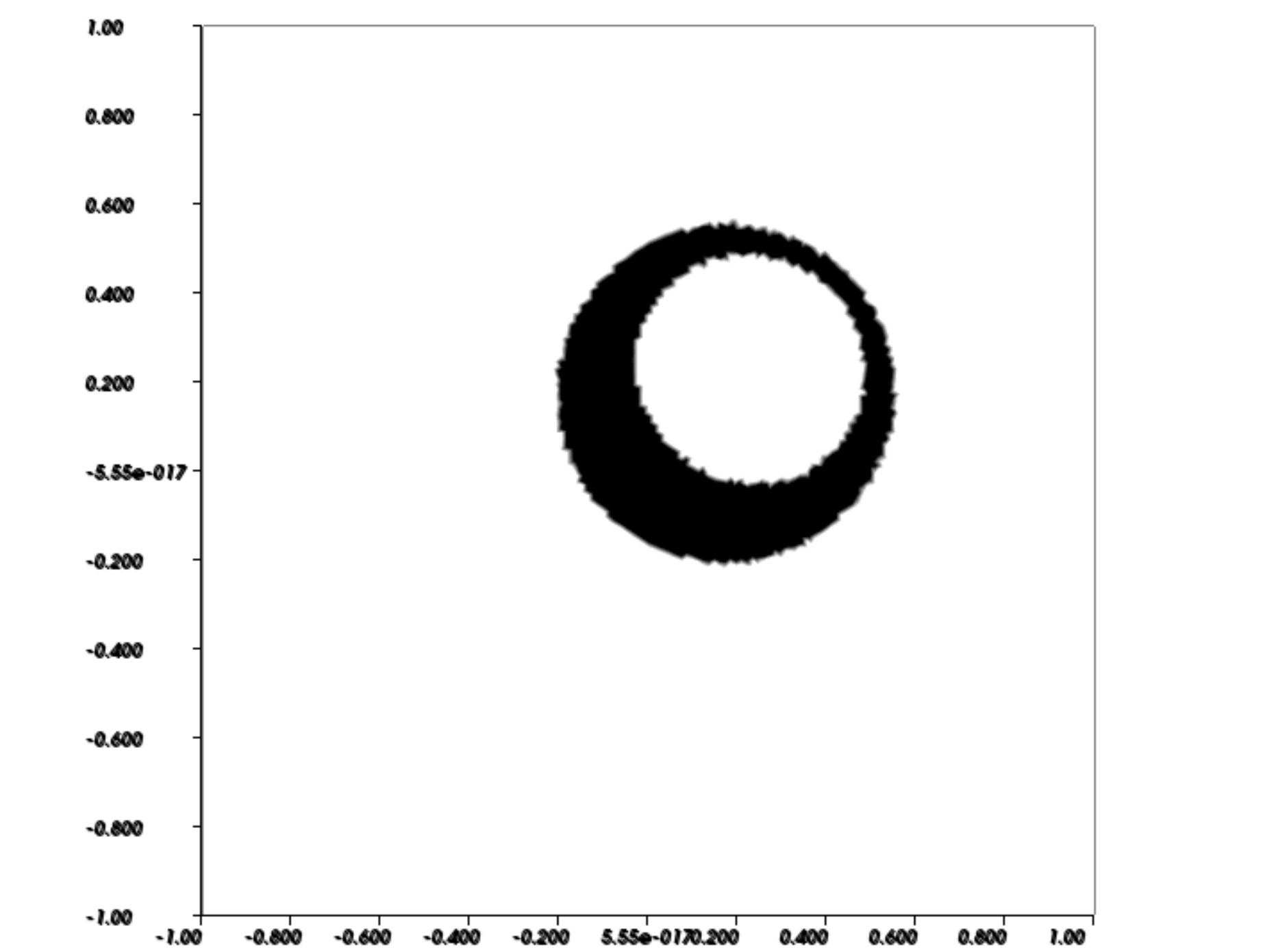}
\
\includegraphics[width=4.5cm]{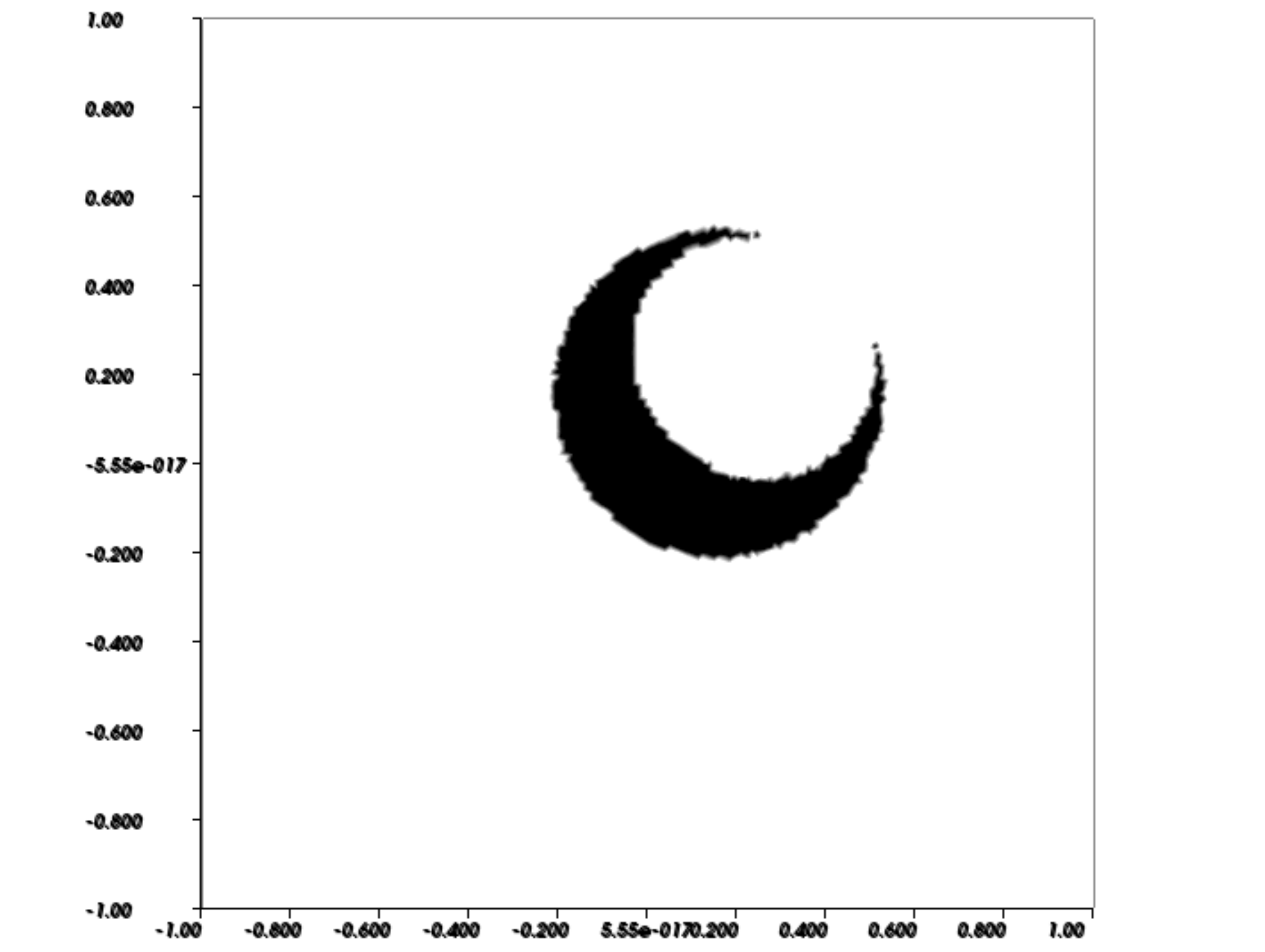}
\
\includegraphics[width=4.5cm]{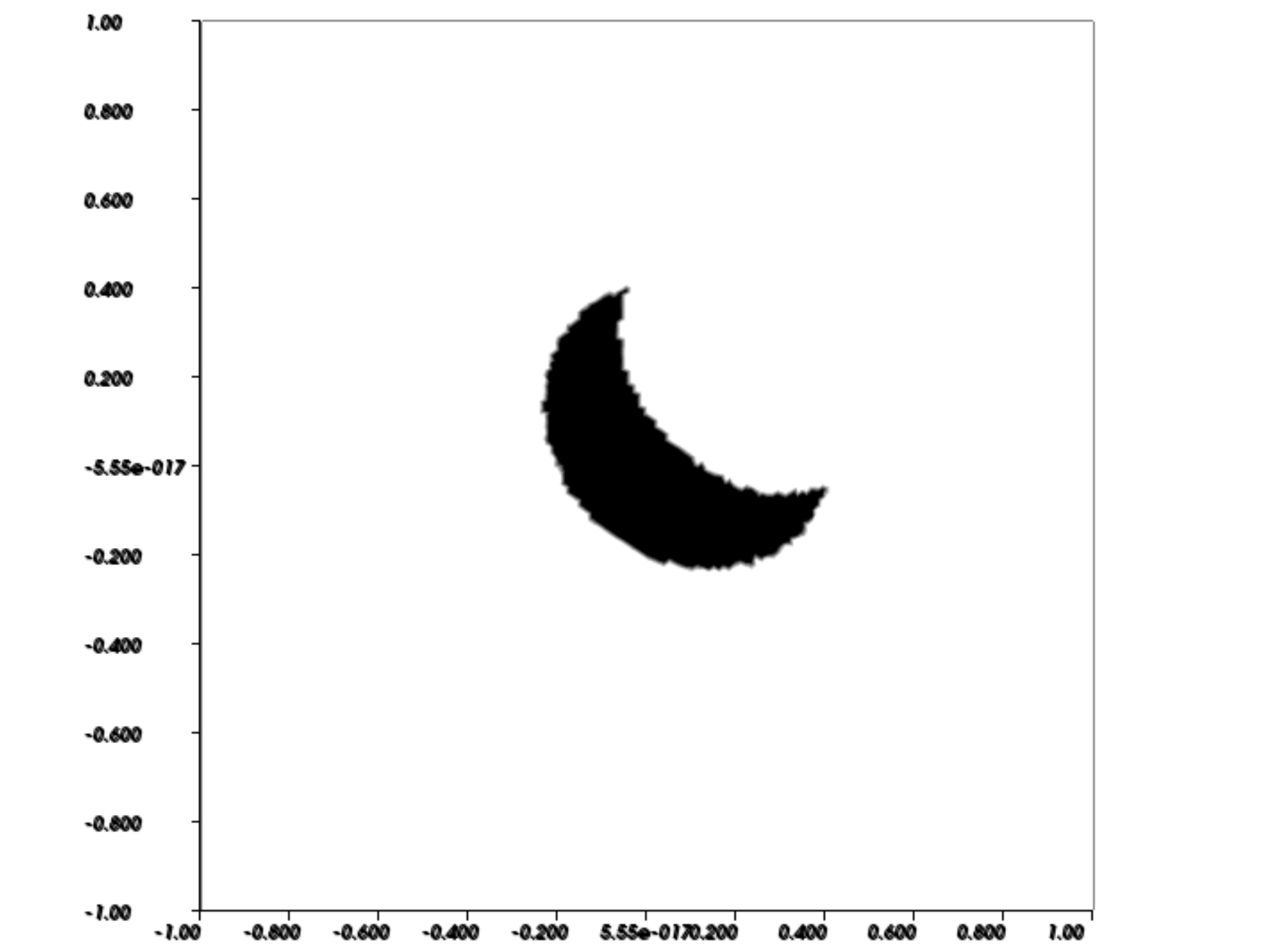}
\end{center}
\caption{Example 1b. Initial domain (top, left), intermediate
and the final domain (bottom, right).
\label{fig:ex2_xi}}
\end{figure}

The stopping test is obtained for $k=5$.
The objective function (\ref{1.3}) is $0.455005$
for the solution of the elliptic system (\ref{1.2})-(\ref{1.3})
in the domain $\Omega_g$ and $0.204318$ for the solution in the domain of boundary
the zero level sets of $y$.
The domain changes its topology,
the initial domain is double connected and the final one is simply connected,
see Figure \ref{fig:ex2_xi}. The penalization term is here a sum of two integrals as explained after (\ref{3.4}).
The corresponding values of the objective function (\ref{3.1})
are reported in Table \ref{tab:ex1b_J}.

\begin{table}[ht]
\begin{center}
\begin{tabular}{|c|c|c|c|c|c|c|}
\hline
iteration & initial & k=2 &   &  &  & final \\  \hline
$t_1$ & 1.57612 & 2.03842 & 2.10354 & 2.12639 & 1.70448 & 0.494585  \\  \hline
$t_2$ & 3.70008 & 0.147847 & 0.123044 & 0.106989 & 0.047237 &0.002589  \\  \hline
$J$   & 37002.4 & 1480.51 & 1232.54 & 1072.02 & 474.078 & 26.3898\\  \hline
\end{tabular}
\end{center}
\caption{Example 1b. The computed objective function (\ref{3.1}), i.e.
$J=t_1+\frac{1}{\epsilon}t_2$, where
  $t_1$ and $t_2$ are as before. The columns 4, 5, 6 corespond
  to intermediate configurations
obtained during the line-search after k=2.}
\label{tab:ex1b_J}
\end{table}

\end{document}